\documentclass[11pt]{amsart}

\usepackage{amsmath}
\usepackage{amsxtra}
\usepackage{amscd}
\usepackage{amsthm}
\usepackage{amsfonts}
\usepackage{amssymb}
\usepackage{mathtools}
\usepackage{eucal}
\usepackage[all]{xy}
\usepackage{graphicx}
\usepackage{tikz}
\usetikzlibrary{shapes.misc}
\usepackage{tikz-cd}
\usepackage{mathrsfs}
\usepackage{euler}
\usepackage{color}
\usepackage{caption}
\usepackage{enumitem}
\usepackage{faktor}
\usepackage{bbm}
\usepackage{stmaryrd}
\usepackage[T1]{fontenc} 
\usepackage[utf8]{inputenc}%
\usetikzlibrary{matrix, arrows, decorations.pathmorphing}
\usepackage[hidelinks]{hyperref}
\usepackage{float}
\usepackage{multicol}
\usepackage{scrextend}
\usepackage{lineno}

\usepackage[backref=true,
                        style=alphabetic,
                        doi=false,
                        maxcitenames=8,
                        maxnames=50,
                        maxalphanames=4]{biblatex}
\usepackage[colorinlistoftodos, textsize=tiny]{todonotes}

\def \vkv{{\begin{xy}(0, 1)*+{v}="A", (10, 1)*+{v'}="B",
\ar@{-}^k"A";"B"
\end{xy}}}

\def \tkt{{\begin{xy}(0, 1)*+{t}="A", (10, 1)*+{t'}="B",
\ar@{-}^k"A";"B"
\end{xy}}}

\theoremstyle{plain}
    \newtheorem{theorem}{Theorem}
    \newtheorem{proposition}[theorem]{Proposition}
    \newtheorem{lemma}[theorem]{Lemma}
    \newtheorem{lemdef}[theorem]{Lemma-Definition}
    \newtheorem{corollary}[theorem]{Corollary}

    \newtheorem{maintheorem}{Theorem}

    \newtheorem{definition}[theorem]{Definition}
    \newtheorem{example}[theorem]{Example}
    \newtheorem{conjecture}[theorem]{Conjecture}

\theoremstyle{remark}
    \newtheorem{remark}[theorem]{Remark}
    
    \newtheorem*{remark*}{Remark}
    \newtheorem{introremark}{Remark}[maintheorem] 

\DeclareMathAlphabet{\mathcal}{OMS}{cmsy}{m}{n}

\newcommand{\A}{\mathbb{A}}

\newcommand{\C}{\mathbb{C}}

\newcommand{\F}{\mathbb{F}}
\newcommand{\G}{\mathbb{G}}

\newcommand{\LL}{\mathbb{L}}

\renewcommand{\P}{\mathbb{P}}
\newcommand{\Q}{\mathbb{Q}}
\newcommand{\R}{\mathbb{R}}
\newcommand{\T}{\mathbb{T}}
\newcommand{\Z}{\mathbb{Z}}
\newcommand{\Zpos}{\mathbb{Z}_{>0}}
\newcommand{\Znn}{\mathbb{Z}_{\geq 0}}

\newcommand{\calb}{\mathcal{B}}

\newcommand{\calf}{\mathcal{F}}
\newcommand{\calg}{\mathcal{G}}

\newcommand{\calm}{\mathcal{M}}
\newcommand{\calo}{\mathcal{O}}

\newcommand{\calu}{\mathcal{U}}

\newcommand{\rmn}{\mathrm{N}}

\newcommand{\Aut}{{\mathrm{Aut}}}

\newcommand{\Id}{\mathrm{Id}}

\newcommand{\Bl}{\operatorname{Bl}}

\newcommand{\Frieze}{\mathrm{Frieze}}

\renewcommand{\div}{{\mathrm{div}}}

\newcommand{\ord}{{\mathrm{ord}}}

\newcommand{\paren}[1]{\mathopen{}\left(#1\right)\mathclose{}}
\newcommand{\set}[1]{\mathopen{}\left\{#1\right\}\mathclose{}}

\newcommand{\verts}[1]{\mathopen{}\left\lvert#1%
\right\rvert\mathclose{}}

\newcommand\restr[2]{{
    \left.\kern-\nulldelimiterspace%
    #1 
    \right|_{#2} 
    }}

\makeatletter
\AtBeginDocument
    {
        \def\@thm#1#2#3{%
            \ifhmode
                \unskip\unskip\par
            \fi
            \normalfont
            \trivlist
            \let\thmheadnl\relax
            \let\thm@swap\@gobble
            \let\thm@indent\indent 
            \thm@headfont{\scshape}
            \thm@notefont{\fontseries\mddefault\upshape}%
            \thm@headpunct{.}
            \thm@headsep 5\p@ plus\p@ minus\p@\relax
            \thm@space@setup
            #1
            \@topsep \thm@preskip 
            \@topsepadd \thm@postskip 
            \def\dth@counter{#2}%
            \ifx\@empty\dth@counter
                \def\@tempa{%
                    \@oparg{\@begintheorem{#3}{}}[]%
                }%
            \else
                \H@refstepcounter{#2}%
                \hyper@makecurrent{#2}%
                \let\Hy@dth@currentHref\@currentHref
                \AddToHookNext{para/begin}{\MakeLinkTarget*{%
\Hy@dth@currentHref}}%
                \def\@tempa{%
                    \@oparg{\@begintheorem{#3}{\csname the#2\endcsname}}[]%
                }%
            \fi
            \@tempa
        }%
    \dth@everypar={%
        \@minipagefalse \global\@newlistfalse
        \@noparitemfalse
        \if@inlabel
            \global\@inlabelfalse
            \begingroup \setbox\z@\lastbox
             \ifvoid\z@ \kern-\itemindent \fi
            \endgroup
            \unhbox\@labels
        \fi
        \if@nobreak \@nobreakfalse \clubpenalty\@M
        \else \clubpenalty\@clubpenalty \everypar{}%
        \fi
    }%
}

\newcommand{\subalign}[1]{%
    \vcenter{%
        \Let@ \restore@math@cr \default@tag
        \baselineskip\fontdimen10 \scriptfont\tw@
        \advance\baselineskip\fontdimen12 \scriptfont\tw@
        \lineskip\thr@@\fontdimen8 \scriptfont\thr@@
        \lineskiplimit\lineskip
        \ialign{\hfil$\m@th\scriptstyle##$&$\m@th\scriptstyle{}%
##$\hfil\crcr
            #1\crcr
        }%
    }%
}

\makeatother

\usepackage{cleveref}
\crefname{maintheorem}{theorem}{theorems}
\Crefname{maintheorem}{Theorem}{Theorems}

\title{Positive integral points on Mordell--Schinzel surfaces}

\author{Antoine de Saint Germain}
\address{
Department of Mathematics and New Cornerstone Science Laboratory \\
The University of Hong Kong \\
Pokfulam Road \\
Hong Kong}
\email{adsg96@hku.hk}

\author{Robin Zhang}
\address{Institut des Hautes \'{E}tudes Scientifiques \\
35 route de Chartres \\
91893 Bures-sur-Yvette \\ France}
\email{rzhang@ihes.fr}

\date{\today}

\begin{document}


\begin{abstract}
    We apply the theory of generalised cluster algebras to classify the finiteness of positive integral points on the class of log K3 surfaces defined by $xyz = G(x, y)$ where $G$ has positive coefficients. We furthermore establish results on Zariski density, number of orbits and asymptotic growth of such points. These provide a Diophantine realisation of the finite–affine–indefinite trichotomy of rank-two root systems.
\end{abstract}

\maketitle

\setcounter{tocdepth}{1}
\tableofcontents


\section{Introduction}
\label{sec:introduction}

\subsection{Mordell--Schinzel varieties}
Let $K$ be a field of characteristic $0$, let $d \geq 1$,
and let $G \in K[x_1, \ldots, x_d]$. We call
\[
    Y_G =
    \set{(x_0, \ldots, x_d) \in \A_K^{d + 1}
        : x_0x_1\cdots x_d = G(x_1, \ldots, x_d)}
\]
the $d$-dimensional Mordell--Schinzel variety defined by $G$. When $G$ has integer coefficients, the same equation defines
an integral model. We write $Y_G(\Z)$ for its integral points
and $Y_G(\Zpos)$ for those with all coordinates positive;
here $\Zpos = \{1, 2, \ldots\}$ and
$\Znn = \{0, 1, 2, \ldots\}$.
On the locus $x_1\cdots x_d \neq 0$,
projection onto the last $d$ coordinates
identifies integral points
with tuples of nonzero integers satisfying
$x_1\cdots x_d \mid G(x_1, \ldots, x_d)$.
This relates $Y_G(\Z)$
to the polynomial divisibility problem
considered by
\cite{schinzel-2000,corvaja-zannier-2010,corvaja-zannier-new}:
for $F, G \in \Z[x_1, \ldots, x_d]$
with $F \neq 0$, when does $G/F$ take integer values at
infinitely many points of $\Z^d$?

Mordell--Schinzel varieties are closely related to several classical
Diophantine equations. For example, the Markoff surface
$3xyz = x^2 + y^2 + z^2$ is the fiber $x_0 = 3$ of the
Mordell--Schinzel threefold defined by
$G = x_1^2 + x_2^2 + x_3^2$.
More generally, fixing $x_0$ gives the Markoff--Hurwitz
varieties studied in
\cite{baragar:Markoff-hurwitz, chang-shparlinski,
    gamburd-asymptotics}, as well as the Rosenberger equations
\cite{rosenberger:equation}.
Fricke cubic surfaces and affine cones over members of the
Dwork family also arise as such fibers; see
\cite{cantat-character-varieties, Whang:log-CY, katz-dwork,
    ghk-automorphy}.

Smooth and geometrically integral Mordell--Schinzel varieties are log Calabi--Yau
in the sense of Gross--Hacking--Keel
\cite[Definition~1.1]{GHK-birational}, as we demonstrate in Proposition~\ref{prop:MS-log-CY}.
In the case of a smooth Mordell--Schinzel surface $Y_G$
with $G(0, 0) \neq 0$,
this means that there exists
a smooth projective surface $X$ and
a reduced simple normal crossings divisor $D$
such that $X \setminus D = Y_G$ and $K_X + D \sim 0$,
where $K_X$ is the canonical divisor of $X$
and $\sim$ denotes linear equivalence.
Even if $Y_G$ is singular, a resolution of $Y_G$
has logarithmic Kodaira dimension $0$.

The logarithmic version of Campana's conjecture predicts
that integral points on log Calabi--Yau varieties are
potentially Zariski dense; see
\cite[Sections~1.2.4 and~1.3.1]{corvaja-zucconi}.
However, there is no general prediction
for the finiteness or Zariski density
of integral points on a given model
of log Calabi--Yau varieties
over a fixed ring of integers.
In this article, we study these questions for
integral points on Mordell--Schinzel surfaces over $\Z$.

\subsection{The Mordell--Schinzel program}
\label{sec:MSproblem}
For surfaces, we use coordinates $(x, y, z)$.
In 1939, Jacobsthal~\cite[\S2, Satz~5]{jacobsthal} showed
that the equation
\[
    xyz = x^2 + x + y^2 + y - 1
\]
has exactly four integer solutions.
In 1952, Mordell~\cite[Theorem~II]{mordell-1952}
developed a recursive construction giving rise to infinitely many integer solutions to
\[
    xyz = ax^3 + by^3 + c
\]
and claimed~\cite[p.~79]{mordell-1952} that his method
applies to arbitrary $G$, although Jacobsthal's example
already contradicts this broader assertion.
Barnes~\cite{barnes}, Mills~\cite{mills2}, and
most notably Schinzel~\cite{schinzel-2015, schinzel-2018} later
re-examined Mordell's construction, completed his cubic argument,
and established the infinitude of integer solutions
to $xyz = G(x, y)$
for a large family of polynomials $G \in \Z[x, y]$.
These results led Koll\'{a}r--Li~\cite[Conjecture~2]{kollar-li}
to formulate a general conjecture, which we now recall.

Write $A(x) := G(x, 0)$ and $B(y) := G(0, y)$,
set
\[
    m := \deg A \quad \text{and} \quad n := \deg B,
\]
and call $(m, n)$ the bidegree of $G$
with the convention that $\deg 0 := 0$.
Let $a$ and $b$ be the coefficients of $x^m$ and $y^n$
in $A$ and $B$ respectively,
and let $c := G(0, 0)$ be the constant term of $G$.

\begin{conjecture}[{Mordell--Schinzel conjecture \cite[Conjecture~2]{kollar-li}}]
    \label{conj:mordell-schinzel}
    Fix $G \in \Z[x, y]$ with bidegree $(m, n)$
    and $c \neq 0$.
    If $m, n \geq 3$, then $Y_G(\Z)$ is infinite.
\end{conjecture}

Koll\'{a}r--Li~\cite[Theorem~3 and Paragraph~8]{kollar-li}
prove that $Y_G(\Z)$ is infinite if $m = n = 3$
and that $Y_G(\Z)$ is Zariski dense
if $m = n = 3$ and $a = b = c = 1$,
using isomorphisms among
a finite family of auxiliary surfaces
and the automorphisms studied by
Koll\'{a}r--Villalobos-Paz~\cite[Theorem~28]{KVP:cubic}.
Koll\'{a}r--Li~\cite[p.4]{kollar-li} expect
Zariski density to hold
for all $m,n \geq 3$.

\subsection{Main results}
\label{ss:main-result}
In this paper, we study positive integral points on $Y_G$
when $G \in \Znn[x,y]$.
For convenience, we assume {\it throughout the article}
that $G$ is not identically zero.
Let $\Aut_\Z(Y_G)$ be the group of automorphisms
of $Y_G$ defined over $\Z$,
and consider the subgroup of positivity-preserving automorphisms
\[
    \Aut_\Z^+(Y_G) :=
    \set{\varphi \in \Aut_\Z(Y_G)
        : \varphi\big(Y_G(\Zpos)\big) = Y_G(\Zpos)}.
\]

\begin{maintheorem}
    \label{thm:trichotomy}
    Fix $G \in \Znn[x, y]$ with bidegree $(m, n)$ and $c \neq 0$.
    \begin{enumerate}[label*=(\roman*)]
        \item If $mn < 4$, then $Y_G(\Zpos)$ is finite.
        \item If $mn = 4$, then $Y_G(\Zpos)$ is infinite. Furthermore, if $abc^{m+n-4} = 1$, then $Y_G(\Zpos)$ is not Zariski dense and consists of only finitely many $\Aut_\Z^+(Y_G)$-orbits.
        \item If $mn > 4$, then $Y_G(\Zpos)$ is Zariski dense.
    \end{enumerate}
\end{maintheorem}
\begin{introremark}
    Our proof of Theorem~\ref{thm:trichotomy}(i)
    is effective in the sense that we obtain an explicit bound on $\#Y_G(\Zpos)$.
    There are specific cases in which
    $\#Y_G(\Zpos)$ is precisely known,
    e.g. for special cases of bidegree $(3, 1)$ in
    \cite{mohanty-1977, ramasamy-mohanty96, subburam-togbe:G2}
    and for the family $Y_{x^m + y + 1}$
    in \cite[Corollary A]{zhang-dynkin}.
\end{introremark}

Theorem~\ref{thm:trichotomy} classifies the
finiteness of $Y_G(\Zpos)$ for $G$ with constant term $c \neq 0$
solely in terms of the bidegree $(m, n)$ of $G$.
Our next result extends this finiteness classification
to $G$ with $c = 0$.
When $Y_G$ is geometrically integral, so that $A$ and $B$
are nonzero, set $k := \ord_0 A$ and $\ell := \ord_0 B$,
where $\ord_0$ denotes the order of vanishing at $0$.

\begin{maintheorem}
    \label{thm:finiteness-classification}
    Fix $G \in \Znn[x, y]$ with bidegree $(m, n)$ and $c = 0$.
    \begin{enumerate}[label*=(\roman*)]
        \item If $Y_G$ is geometrically integral,
              then $Y_G(\Zpos)$ is finite if and only if
              $\min(m, n) = 1$ and
              $\max(m, n) - \max(k, \ell) < 4$.
        \item Otherwise, $Y_G(\Zpos)$ is infinite.
    \end{enumerate}
\end{maintheorem} 

Theorem \ref{thm:trichotomy} also reveals an $mn$-trichotomy,
reflecting the Lie-theoretic trichotomy of
finite, affine, and indefinite
type\footnote{Also called finite, tame and wild type, respectively} root systems;
this connection
is made apparent by the theory of generalised cluster algebras
and is explained in Section \ref{subsec:outline}. 
The $mn = 4$ case corresponds to affine type
and is of particular interest.
Under the assumption that $abc^{m+n-4}=1$, 
our proof of Theorem~\ref{thm:trichotomy}(ii)
yields a point-counting asymptotic.
We measure the size of a positive integral point
by the height function
$h(x, y, z) = \max(x, y)$ and, for $N \geq 1$, set
\[
    E_G(N) = \#\set{P \in Y_G(\Zpos) : h(P) \leq N}.
\]

\begin{maintheorem}
    \label{thm:affine-counting}
    Fix $G \in \Znn[x, y]$ with bidegree $(m, n)$ and $c \neq 0$. If $mn = 4$ and $abc^{m + n - 4} = 1$, then
    there is a constant $C_G > 0$ such that
    $E_G(N) = C_G\log N + O_G(1)$.
\end{maintheorem}

For the remaining affine type cases ($mn = 4$ and $abc^{m+n-4} > 1$)
and for indefinite types ($mn > 4$),
we propose the following conjecture.

\begin{conjecture}
    \label{conj:indefinite-orbits}
    Fix $G \in \Znn[x, y]$ with bidegree $(m, n)$ and $c \neq 0$.

\begin{enumerate}[label*=(\roman*)]
    \item If $mn = 4$ and $abc^{m+n-4} > 1$, then $Y_G(\Zpos)$ is Zariski dense and consists of infinitely many $\Aut_\Z^+(Y_G)$-orbits.
    \item If $mn > 4$, then $Y_G(\Zpos)$
    consists of infinitely many $\Aut_\Z^+(Y_G)$-orbits.
\end{enumerate}
\end{conjecture}

We prove the density assertion of Conjecture \ref{conj:indefinite-orbits}(i)
for $G(x, y) = ax^2 + by^2 + c$ with $ab > 1$
in Proposition \ref{prop:nonmonic-quadratic-density},
and prove Conjecture \ref{conj:indefinite-orbits}(ii)
for $G(x, y) = x^3 + y^2 + 1$ in  Proposition~\ref{prop:indefinite-3-2-orbits}.

\subsection{Proof overview and organization}
\label{subsec:outline}
In Section~\ref{sec:generalities}, we identify the open subset
$x_1 \cdots x_d \neq 0$ of $Y_G$ with the torus $\G_m^d$ of
nonzero coordinate tuples.
For smooth geometrically integral $Y_G$,
we extend the differential form
$(dx_1 \wedge \cdots \wedge dx_d)/(x_1\cdots x_d)$
to a nowhere-vanishing regular $d$-form on $Y_G$.
The embedding of the torus $\G_m^d$ and
this extension of the differential form
imply the log Calabi--Yau property by
Gross--Hacking--Keel~\cite[Lemma~1.4(i)]{GHK-birational}
(Proposition~\ref{prop:MS-log-CY}).

We now turn to our arithmetic results on Mordell--Schinzel surfaces. For $G \in \Znn[x,y]$ with $m,n \geq 1$ and $c \neq 0$, there is a unique $H \in \Znn[x, y]$ such that
\[
    G(x, y) = A(x) + B(y) - c + xyH(x, y).
\]
We call both $G$ and $Y_G$ pure when $H = 0$, and write $Y_G = Y_{A,B}$ in this case. The substitution $z \mapsto z - H(x, y)$ identifies $Y_G$ with the pure Mordell--Schinzel surface obtained by setting $H=0$. Moreover, this isomorphism preserves the bidegree,
preserves the height $h$, and induces a bijection
on positive integral points (Lemma \ref{lem:mixed-reduction}). Thus, we may assume $G$ is pure.

The space of pure Mordell--Schinzel surfaces of bidegree $(m,n)$ is naturally parametrised by (an open subset of) $\calb_{m,n}$, the affine space of coefficients of $A$ and $B$ 
(including their common constant term $c$). In Section~\ref{sec:gca}, we work with the flat family
$\pi: \calm_{m,n} \longrightarrow \calb_{m,n}$ with 
$\calm_{m,n}$ defined by
$xyz = A(x) + B(y) - c$, where now $A(t),B(t) \in \Z[\calb_{m,n}][t]$. We prove that $\Z[\calm_{m,n}]$ is isomorphic to a rank-$2$
generalised upper cluster algebra over $\Z[\calb_{m,n}]$ denoted $\calu_{m,n}$. This equips $\calm_{m,n}$ with a set of regular functions called cluster variables, defined recursively via the so-called exchange relations (see Appendix \ref{sec:appendix}). 

Every (generalised upper) cluster algebra comes with a family of distinguished morphisms called friezes. In our case, roughly speaking, these are $\Z[\calb_{m,n}]$-algebra homomorphisms $\calu_{m,n} \to \Z$ sending every cluster variable to a positive integer. 
Such a frieze is said to be adapted to $G \in \Znn[x, y]$ if it specializes the coordinates of $\calb_{m,n}$
to the coefficients of $G$. In Theorem \ref{thm:frieze-and-integer-points}, we identify 
$Y_G(\Zpos)$ with the set of adapted friezes. This 
translates the finiteness statements of Theorem \ref{thm:trichotomy} into corresponding statements about adapted friezes of $\calu_{m,n}$.


Recall that a (generalised upper) cluster algebra
is said to be of finite cluster-type if it admits only finitely many cluster variables. According to \cite{fomin-zelevinsky-2,chekhov-shapiro}, such algebras are classified by finite root systems. In particular, $\mathcal{U}_{m,n}$ is of finite cluster-type if and only if $mn \leq 3$; the cases $mn=1, 2$ and $3$ correspond to types $A_2,B_2$ and $G_2$, respectively. This naturally leads to the  finiteness classification in
Theorem \ref{thm:trichotomy}. Indeed, when $mn \geq 4$, we use the infinitude of cluster variables to construct infinitely many friezes (called unitary friezes) adapted to each $G$. The existence of these friezes is a consequence of the remarkable Laurent phenomenon. When $mn \leq 3$, we use the 
finiteness of cluster variables to construct regular functions on $Y_G$ whose values on positive integral points are bounded. Together with elementary divisibility arguments when $mn = 0$, this proves the finiteness classification of Theorem \ref{thm:trichotomy}.

To prove Theorem~\ref{thm:trichotomy}(iii), we again use the unitary friezes to
construct distinct positive integral points
$P_k = (x_k, y_k, z_k)$ such that $x_k \rightarrow \infty$ and
$\log y_k / \log x_k$ converges to an irrational number
as $k \rightarrow \infty$.
But if an algebraic curve in $Y_G$
contains infinitely many $P_k$,
then this limit,
after passing to a subsequence,
equals a ratio of the orders of poles
along one branch at infinity, hence is rational, which is a contradiction.

Our proof of Theorem~\ref{thm:finiteness-classification} is more elementary, and relies on explicit constructions and
reductions to Mordell--Schinzel surfaces for which Theorem~\ref{thm:trichotomy} applies. 

The (generalised upper) cluster algebras $\calu_{m,n}$ with $mn \geq 4$ are classified by infinite rank-2 root systems, and are thus split into the affine ($mn = 4$) and indefinite ($mn >4$) types. 
Suppose that $mn = 4$ and $abc^{m+n-4} = 1$,
as in the statement of Theorem~\ref{thm:trichotomy}(ii)
and Theorem~\ref{thm:affine-counting}.
Following Koll\'{a}r--Villalobos-Paz~\cite{KVP:cubic},
we show in Section~\ref{sec:companion} that the rational
maps on $\A^2$ given by 
\begin{equation}\label{eq:one-step-mutation-rational}
    (x,y) \mapsto \left(\frac{B(y)}{x},y\right) \quad \text{and} \quad (x,y) \mapsto \left(x,\frac{A(x)}{y}\right),
\end{equation}
lift to isomorphisms $\mu_{1}: Y_{A,B} \to Y_{A^\dagger,B}$ and $\mu_{2}: Y_{A,B} \to Y_{A,B^\dagger}$ over $\Z$,
with $A^\dagger$ and $B^\dagger$ defined in
Section~\ref{subsec:two-isoms}.
The vertices of the diagram
\[
\begin{tikzcd}
    Y_{A,B} \arrow[d, "\mu_1"] \arrow[r,"\mu_2"] & Y_{A,B^\dagger} \arrow[d, "\mu_1^\dagger"] \\
    Y_{A^\dagger,B} \arrow[r,"\mu_2^\dagger"] & Y_{A^\dagger, B^\dagger}
\end{tikzcd}
\]
are the {\it companion surfaces} of $Y_{A,B}$,
as studied in \cite{KVP:cubic, kollar-li}.
They constitute the objects of a groupoid $\calg_{A,B}$,
whose morphisms are generated by the edges of the diagram above.
We write $\Gamma_{A,B}$ for the vertex group of $Y_{A,B}$, that is, the group of morphisms of $\calg_{A,B}$
that are automorphisms of $Y_{A,B}$.

In Section~\ref{sec:orbits}, we make crucial use of this vertex group to study the number of $\Aut_\Z^+(Y_G)$-orbits
in $Y_G(\Zpos)$. First, we show that each $\Gamma_{A,B}$-orbit has a ``minimal'' element, and we bound its coordinates. This implies there are only finitely many
$\Gamma_{A,B}$-orbits in $Y_G(\Zpos)$. Since
$\Gamma_{A,B} \subseteq \Aut_\Z^+(Y_G)$, there are also
only finitely many $\Aut_\Z^+(Y_G)$-orbits.
Moreover, each $\Gamma_{A,B}$-orbit of $Y_G(\Zpos)$ lies on a level
curve of an explicit nonconstant regular function of $Y_G$, as predicted by the scattering diagram of $\calu_{m,n}$ defined in \cite{GHKK,burcroff-Lee-mou:generalised}. Labelling the 
points in each $\Gamma_{A,B}$-orbit suitably reveals that their coordinates satisfy linear recurrence relations, and thus the orbit has counting asymptotic $C\log N + O(1)$, where $C$ is an explicit positive
constant depending on the orbit.
Together, these prove Theorem~\ref{thm:trichotomy}(ii)
and Theorem \ref{thm:affine-counting}.

Section~\ref{sec:examples} gives examples
and applications to binary quadratic forms.
Appendix~\ref{sec:appendix} records the definitions and
properties of generalised cluster algebras used in the proofs.

\subsection*{Acknowledgments}
The authors thank Pietro Corvaja, Benjamin Kane, Jennifer Li,
Peter Sarnak, Katherine Stange, Junho Whang, and Umberto Zannier for
helpful discussions.
A. de Saint Germain was partially supported by
the New Cornerstone Science Foundation
through the New Cornerstone
Investigator Program awarded to Professor Xuhua He.
R. Zhang was partially supported
by the National Science Foundation under Grant No.~DMS-2303280.

\subsection*{AI usage}
Most of this article was
close to its current form in early 2026,
when the results were first announced,
and the results were developed without the use of AI.
Since then, the authors have used ChatGPT 5.6 and 6
to help with the final revisions,
harmonisation of notation,
and error-verification.

\numberwithin{equation}{section}
\numberwithin{theorem}{section}

\section{The geometry of Mordell--Schinzel varieties}
\label{sec:generalities}
We study the geometry of $Y_G$ and prove that it is log
Calabi--Yau when it is smooth and geometrically integral.
Throughout this section, $K$ has characteristic zero, $d \geq 1$, $Q = x_1\cdots x_d$, $G \in K[x_1, \ldots, x_d]$ and $F = x_0Q - G$. We write $V(I)$ for the closed subscheme defined by an ideal $I$, and $D(f)$ for the principal open set where $f \neq 0$. The geometric constructions below are variations of standard constructions in \cite{GHK-birational,blanc-dolgachev,KVP:cubic}; we record the statements needed later.

\subsection{Affine modifications}
\label{sec:affine-mod}
Since $x_0Q-G$ has degree one in $x_0$, Gauss's lemma gives
\begin{equation}\label{eq:integrality-criterion}
  \begin{split}
  Y_G\text{ is geometrically integral}
  &\quad\Longleftrightarrow\quad \gcd(Q,G)=1\\
  &\quad\Longleftrightarrow\quad x_i\nmid G\quad(1\leq i\leq d).
  \end{split}
\end{equation}
Henceforth, assume $Y_G$ is geometrically integral.

Projection onto $(x_1, \ldots, x_d)$ identifies the open set $D(Q)\subset Y_G$ with the torus $(\G_m)^d$. In particular, $Y_G$ is rational, and every singular point lies in $V(Q)$. 
In $V(Q)$, the defining equation of $Y_G$ implies $G = 0$. In other words, there is a natural isomorphism $Y_G \cap V(Q) \cong \A^1_{x_0} \times V(Q, G)$. This suggests a more geometric construction of $Y_G$, as follows. 

\begin{proposition}
    \label{prop:blowup-general}
    Assume that $\gcd(Q, G) = 1$. 
    If $Z = \Bl_{(Q, G)}\A_K^d$ and $\widetilde D$ is the
    strict transform of $D = V(Q)$, then
    $Y_G \cong Z \setminus \widetilde D$. In particular, 
    $K[Y_G] = R[G/Q] \subset R[1/Q]$ where $R = K[x_1,\ldots, x_d]$. 
\end{proposition}

\begin{proof}
    Since $K[Y_G]$ is a domain, allowing division by $Q$ embeds
    it into $R[1/Q]$. The relation $x_0 = G/Q$ identifies
    its image with $R[G/Q]$.

    The blow-up has equation $QU = GV$ in
    $\A_K^d \times \P^1_{[U : V]}$, as in
    Kaliman--Zaidenberg
 ~\cite[Example~1.4]{Kaliman-Zaidenberg:affine-mod}.
    Its intersection with $V = 0$ is a copy of $D$.
    Since no component of $D$ is contained in $V(Q, G)$,
    this intersection is $\widetilde D$.
    In the complementary chart $V \neq 0$, setting
    $x_0 = U/V$ gives $Qx_0 = G$.
\end{proof}
In the terminology of \cite{Kaliman-Zaidenberg:affine-mod}, $Y_G$ arises as an {\it affine modification of $\A^d$ along} the divisor $V(Q)$ with {\it center} $V(Q,G)$. Variations of this construction have appeared in \cite{blanc-dolgachev,GHK-birational,KVP:cubic}. When $d = 2$, the center $V(xy,G)$ is supported at only finitely many points of the coordinate axes. 

\begin{remark}
For every $H \in R$, the map
$(x_0, x_1, \ldots, x_d) \mapsto (x_0 + H, x_1, \ldots, x_d)$ identifies $Y_G$ with $Y_{G + QH}$. This is clear from the affine modification construction, since $G$ and $G + QH$ coincide as functions on $V(Q)$. 
\end{remark}

\subsection{The log volume form}
For each $1 \leq i \leq d$, write
$G_i = \partial G/\partial x_i$ and
$Q_i = \prod_{j \neq i}x_j$. 

\begin{lemma}
    \label{lem:singular-locus-general}
    A geometric point of $Y_G$ is singular if and only if
    $Q = G = 0$ and $x_0Q_i = G_i$ for every $1 \leq i \leq d$.
\end{lemma}

\begin{proof}
    The partial derivatives of $F$ are $Q$ and $x_0Q_i - G_i$.
    When $Q = 0$, the equation $F = 0$ also gives $G = 0$.
\end{proof}

Recall that an affine integral variety is \emph{normal} when its coordinate
ring is integrally closed, and that a \emph{prime divisor} is an irreducible subvariety of codimension
one.

\begin{lemma}
    \label{lemma:normality-general}
    Assume that $Y_G$ is geometrically integral. Then $Y_G$
    is normal if and only if $G \notin (x_i, x_j)^2$ for
    every $1 \leq i < j \leq d$.
    In particular, $Y_G$ is normal whenever $G(0,\ldots , 0) \neq 0$. 
\end{lemma}

\begin{proof}
    By Serre's criterion for hypersurfaces
    (see~\cite[Tags~00SB and~031S]{stacks-project}),
    it suffices to
    check regularity at the generic points of prime divisors.
    Where exactly $r$ of $x_1, \ldots, x_d$ vanish, the
    equation is $G = 0$ and $x_0$ is free, so the dimension
    is at most $d - r + 1$.
    Thus three or more vanishing coordinates give no divisor.

    If only $x_i$ vanishes, $Q_i$ is nonzero and $x_0$ is
    an independent coordinate on each boundary component.
    Hence $x_0Q_i - G_i$ is nonzero at its generic point.

    If exactly $x_i$ and $x_j$ vanish, the locus is a divisor
    only when $G \in (x_i, x_j)$. In the two transverse
    variables $x_i, x_j$, the term $x_0Q$ has order at least
    two. Thus the generic point is regular precisely when
    $G$ has a nonzero linear term in these variables, namely
    when $G \notin (x_i, x_j)^2$. Its coefficient remains
    nonzero in the fraction field of the other coordinates.
    Finally, a nonzero constant term of $G$ excludes every
    ideal $(x_i, x_j)^2$.
\end{proof}

The canonical bundle of a smooth $d$-dimensional
variety is the line bundle of regular $d$-forms.
For a normal variety, the zeros and poles of a nonzero rational
$d$-form along prime divisors define a \emph{canonical divisor}.
Thus a form with no zeros or poles gives the zero canonical
divisor. 

\begin{proposition}
    \label{prop:canonical-form}
    The form
    \[
        \omega_G =
        \frac{dx_1 \wedge \cdots \wedge dx_d}{x_1\cdots x_d}
        = d\log x_1 \wedge \cdots \wedge d\log x_d
    \]
    on $Q \neq 0$ extends to a nowhere-vanishing regular
    $d$-form on the smooth locus of $Y_G$.
    If $Y_G$ is smooth, its canonical bundle is trivial.
    If $Y_G$ is normal, then $\div_{Y_G}(\omega_G) = 0$,
    so its canonical divisor can be chosen to be $0$.
\end{proposition}

\begin{proof}
    On the chart $\partial F/\partial x_i \neq 0$, the relation
    $dF = 0$ expresses $dx_i$ in terms of the other differentials.
    Therefore
    \[
        \frac{(-1)^i dx_0 \wedge \cdots \wedge
        \widehat{dx_i} \wedge \cdots \wedge dx_d}
        {\partial F/\partial x_i}
    \]
    is a regular nonzero local generator of the $d$-forms;
    the hat $\widehat{dx_i}$ means that $dx_i$ is omitted. The same relation
    $dF = 0$ shows that these expressions agree on overlaps.
    This is the residue calculation in hypersurface adjunction;
    see~\cite[Tags~0B4B and~0AA4]{stacks-project}.
    Since $\partial F/\partial x_0 = Q$, the expression for
    $i = 0$ is $\omega_G$.
    If $Y_G$ is normal, the proof of
    Lemma \ref{lemma:normality-general} shows that the
    generic point of every divisor is smooth. Hence the
    divisor $\div_{Y_G}(\omega_G)$ is zero.
\end{proof}

A smooth compactification of a smooth
variety $U$ is a smooth projective variety $X$ containing
$U$ as a dense open subset. We choose one for which the
boundary $B = X \setminus U$ is a divisor with \emph{simple
normal crossings}.
A logarithmic top form is regular on $U$ and has poles of
order at most one along $B$. More generally, a
\emph{logarithmic pluricanonical form of order $r$} is a
section of the $r$th tensor power of the canonical bundle,
regular on $U$ and with poles of order at most $r$ along $B$. 

Suppose that $U$ is smooth, quasi-projective and geometrically integral. We follow Gross--Hacking--Keel
\cite[Definition~1.1]{GHK-birational}, and call $U$
\emph{log Calabi--Yau} if it has a nowhere-vanishing top form
$\omega$ such that, for every $r \geq 1$, the logarithmic
pluricanonical forms of order $r$ form the space
$K\omega^{\otimes r}$.

\begin{proposition}
    \label{prop:MS-log-CY}
    Suppose that $Y_G$ is geometrically integral.
    \begin{enumerate}
        \item If $Y_G$ is smooth, it is log Calabi--Yau.
        \item If $G(0,\ldots,0)\neq0$, then $Y_G$ has at worst canonical
        singularities, and a resolution has logarithmic Kodaira dimension zero.
    \end{enumerate}
\end{proposition}
\begin{proof}
Both assertions use the open subset $T := D(Q)$ of $Y_G$.
	Projection onto $(x_1, \ldots, x_d)$ identifies $T$ with
	$\G_m^d$, since $x_0 = G/Q$ there. The restriction of
	$\omega_G$ is $dx_1/x_1 \wedge \cdots \wedge dx_d/x_d$.
	On the compactification $(\P^1)^d$, this form has simple
	poles along the coordinate boundary and no zeros.
	Thus $T$ is log Calabi--Yau: for every $r \geq 1$, its
	logarithmic pluricanonical forms of order $r$ form the
	one-dimensional space $K\omega_G^{\otimes r}$.

	Suppose first that $Y_G$ is smooth.
	By \cite[Lemma~1.4(i)]{GHK-birational}, a smooth variety
	containing a log Calabi--Yau open subset is itself
	log Calabi--Yau if the logarithmic top form on that
	open subset extends to a nowhere-vanishing regular form.
	Proposition~\ref{prop:canonical-form} gives precisely
	this extension of $\omega_G$ to $Y_G$, proving (1).

	For (2), assume that $G(0, \ldots, 0) \neq 0$.
	By Lemma~\ref{lemma:normality-general}, $Y_G$ is normal,
	and Proposition~\ref{prop:canonical-form} allows us to
	choose $K_{Y_G} = \div_{Y_G}(\omega_G) = 0$.
	Let $\rho: \widetilde Y_G \to Y_G$ be a resolution that
	is an isomorphism over the smooth locus. In particular,
	$\widetilde Y_G$ still contains $T$ as an open subset.
	We first prove that $\rho^*\omega_G$ is regular.
	Its orders along exceptional divisors are the coefficients
	of $K_{\widetilde Y_G} - \rho^*K_{Y_G}$, called discrepancies;
	their nonnegativity is the condition for canonical
	singularities. The form $\rho^*\omega_G$ is already regular
	over the smooth locus of $Y_G$, so it suffices to exclude
	poles along exceptional prime divisors.

	Suppose that $\rho^*\omega_G$ has a pole along an
	exceptional prime divisor $E$. Then
	\cite[Lemma~1.8]{GHK-birational} implies that the valuation
	$v := \ord_E$ is nonzero and monomial in
	$x_1, \ldots, x_d$: the value of a nonzero polynomial
	is the minimum of the values of its nonzero monomials.
	Since the pullbacks of $x_0, \ldots, x_d$ are regular,
	$v(x_i) \geq 0$ for every $i$. Moreover, some
	$v(x_i)$ with $1 \leq i \leq d$ is positive, since
	otherwise the monomial valuation would be trivial.
	The nonzero constant term of $G$ therefore gives
	\[
		v(G) = 0, \qquad
		v(x_0Q) = v(x_0) + \sum_{i = 1}^d v(x_i) > 0,
	\]
	contradicting $x_0Q = G$.
	Thus $\rho^*\omega_G$ has no poles on $\widetilde Y_G$,
	and $Y_G$ has canonical singularities.

	It remains to compute the logarithmic Kodaira dimension
	of $\widetilde Y_G$. Choose a smooth projective
	compactification with simple normal crossings boundary.
	At the generic point of a boundary divisor, write
	$x_i = t^{e_i}u_i$, where $t$ is a local parameter,
	$e_i \in \Z$, and $u_i$ is a unit. Then
	$dx_i/x_i = e_i\,dt/t + du_i/u_i$.
	In the wedge product defining $\rho^*\omega_G$, terms
	containing $dt/t$ more than once vanish.
	Hence $\rho^*\omega_G$ has at most a simple pole along
	each boundary divisor, and $(\rho^*\omega_G)^{\otimes r}$
	is a nonzero logarithmic pluricanonical form for every
	$r \geq 1$.

	Conversely, restriction to the open torus $T$ injects
	the space of logarithmic pluricanonical forms of order $r$
	on $\widetilde Y_G$ into $K\omega_G^{\otimes r}$.
	Each such space is therefore one-dimensional.
	This proves that $\widetilde Y_G$ has logarithmic
	Kodaira dimension $0$, completing (2).
\end{proof}

\subsection{Positivity}

\begin{lemma}\label{lem:positive-real-component}
    Suppose that $G \in \mathbb{R}_{\geq 0}[x_1,\ldots, x_d]$ and $G(0,\ldots, 0) > 0$. Then 
    $Y_G(\R_{>0}) := Y_G(\R) \cap (\R_{>0})^{d+1}$ is a connected component
    of $Y_G(\R)$. If $G$ has integral coefficients, $Y_G(\Zpos)$ is precisely the set of integral points on this component.
\end{lemma}
\begin{proof}
    Projection to $(x_1, \ldots, x_d)$ identifies $Y_G(\R_{>0})$ with
    $(\R_{>0})^d$, so $Y_G(\R_{>0})$ is connected.
    It is also open in $Y_G(\R)$.

    To show that it is closed, let $(x_0, \ldots, x_d)$ be a limit
    of points in $Y_G(\R_{>0})$. Then $x_i\geq 0$ for $0 \leq i \leq d$.
    Since $G$ has nonnegative coefficients and positive
    constant term, $G(x_1,\ldots, x_d) \geq G(0,\ldots, 0) > 0$.
    The equation $x_0Q = G(x_1, \ldots , x_d)$ therefore forces
    $x_0Q > 0$, so $x_i> 0$ for $0 \leq i \leq d$. Thus $Y_G(\R_{>0})$ is closed as well.
    A nonempty connected subset that is both open and closed
    is a connected component.
\end{proof}

For a totally real number field $L$, let $\calo_L$ be its ring
of integers and set
\[
    \calo_L^+ =\{a \in \calo_L :\sigma(a) > 0 \text{ for every embedding } \sigma: L \hookrightarrow \R\}.
\]
When $G \in \Znn[x_1,\ldots, x_d]$, set $Y_G(\calo_L^+) = Y_G(\calo_L) \cap (\calo_L^+)^{d+1}$. For instance, $Y_G(\calo_L^+) = Y_G(\Zpos)$ when $L = \Q$. 

\begin{example}\label{ex:totally-positive}
    Let $L = \Q(\sqrt{2})$ and $\varepsilon = 1 + \sqrt{2}$. For every nonzero $G \in \Znn[x,y]$ and $n \geq 0$, 
    \[
    \left(\varepsilon^{2n}, \varepsilon^{2n}, \frac{G(\varepsilon^{2n},\varepsilon^{2n})}{\varepsilon^{4n}}\right) \in Y_G(\calo_L^+)
    \]
    Indeed, $\varepsilon^2 = 3 + 2\sqrt{2} \in \calo_L^+$, so each $\varepsilon^{4n}$ is a totally positive unit in $\calo_L$. In particular, $Y_G(\calo_L^+)$ is infinite. Taking $G = 1+ x + y$ shows that Theorem \ref{thm:trichotomy}(i) does not extend to number fields. It would be interesting to study this further.
\end{example}

\subsection{Mordell--Schinzel surfaces}\label{subsec:MS-surface}
Assume now that $d=2$. Fix $G\in K[x,y]$, decomposed as
\[
    G(x, y) = A(x) + B(y) - c + xyH(x, y).
\]
When $c\ne0$, the center
of the affine modification is
\[
 C=V(xy,A(x)+B(y)-c)\subset\A^2.
\]
It is finite and disjoint from the origin. Its points on $y=0$
and $x=0$ correspond to roots of $A$ and $B$, respectively.
By Lemma~\ref{lem:singular-locus-general}, $Y_G$ is smooth if and
only if $A$ and $B$ have only simple roots over $\overline K$.
Under these hypotheses, set $X=\operatorname{Bl}_C\P^2$ and let
$D$ be the strict transform of the coordinate triangle. Then
\begin{equation}\label{eq:triangle-compactification}
 Y_G\cong X\setminus D,\qquad K_X+D\sim0.
\end{equation}

\begin{remark}\label{rem:bidegree}
In the smooth case with $c\ne0$, one may instead use
$\P^1_{[X:U]}\times\P^1_{[Y:V]}$.
The substitution $z\mapsto z-H(x,y)$ identifies $Y_G$ with
$Y_{A(x)+B(y)-c}$ (cf. Lemma~\ref{lem:mixed-reduction}).
The bihomogenisation of this pure polynomial is
\[
 g_h=U^mV^n\bigl(A(X/U)+B(Y/V)-c\bigr)
 \in H^0(\P^1\times\P^1,\mathcal O(m,n)).
\]
Blowing up $V(XY,g_h)$ and removing the strict transform of the
four coordinate divisors gives a surface isomorphic to $Y_G$.

\end{remark}

For singular surfaces, logarithmic Kodaira dimension,
denoted $\overline\kappa$, is computed on a resolution. \begin{lemma}\label{lem:log-kodaira}
Let $Y_G$ be geometrically integral, and let $\widetilde Y_G\to Y_G$
be a resolution that is an isomorphism over its smooth locus.
Then $\overline\kappa(\widetilde Y_G)\leq0$, with equality if $c\neq 0$.
\end{lemma}
\begin{proof}
The resolution contains the torus $D(xy)$; restriction bounds every
logarithmic plurigenus by one \cite[Chapter~11]{iitaka}.
Equality for nonzero constant term is Proposition~\ref{prop:MS-log-CY}(2).
\end{proof}
Thus the finite and infinite cases of
Theorem \ref{thm:trichotomy} both occur at logarithmic
Kodaira dimension zero.

We make the following statement explicit, for future referencing. Its proof is straightforward.

\begin{lemma}
    \label{lem:mixed-reduction}
    Let $G, H \in \Z[x, y]$ and set $J = G + xyH$. The map
    \[
        \begin{tikzcd}[row sep = tiny]
            Y_G \arrow[r, "\phi_H"] & Y_J \\
            (x, y, z) \arrow[r, mapsto] &
            (x, y, z + H(x, y))
        \end{tikzcd}
    \]
    is an isomorphism over $\Z$, which leaves invariant the height $h(x,y,z) = \max(x,y)$. If $G \in \Znn[x,y]$ is nonzero and $H \in \Znn[x,y]$,
    the map $\phi_H$ restricts to a bijection $Y_G(\Zpos) \simeq Y_J(\Zpos)$.
\end{lemma}

\subsection{The lower bound model}
Fix a pure $G \in \Znn[x, y]$, put
\begin{equation}\label{eq:pure-poly}
    G = A(x) + B(y) - c, \qquad A(0) = B(0) = c,
\end{equation}
and write $Y_{A,B} = Y_G$. In this subsection, we think of $Y_G$ as a variety
{\it over} $\Q$. Assume $c \neq 0$ and $m, n \geq 1$. Define
\begin{equation}
    L_G = \Z[u,u',v,v']/I,
\end{equation}
where $I = (uu'-B(v), vv' - A(u))$. We call ${\rm Spec}(L_G)$ the {\it lower bound model}, since it is a straightforward generalisation in rank $2$ of the lower bound model introduced in \cite{fomin-zelevinsky-3} for cluster algebras. When $G = 1 + x^m + y^n$, this was studied by Blanc and Dolgachev \cite{blanc-dolgachev}.
Define the polynomials $A^*,B^* \in \Znn[t]$ by
\begin{equation}\label{eq:A-B-*}
    A(x) -c = xA^*(x) \quad  \text{and} \quad  B(y)-c = yB^*(y).
\end{equation}
There is a ring homomorphism 
\begin{equation}\label{eq:two-form-relation}
\Phi : L_G \longrightarrow \Z[Y_G], \quad (u,u',v,v') \mapsto (x, yz - A^*(x), y, xz-B^*(y)),
\end{equation}
where $\Z[Y_G] = \Z[x,y,z] / (xyz - G)$.

\begin{lemma}\label{lem:2-integral-forms}
    \begin{enumerate}
        \item The extension of scalars along $\Z \to \Z[1/c]$ of the map $\Phi$ given in \eqref{eq:two-form-relation} induces an isomorphism $L_G[1/c] \cong \Z[Y_G][1/c]$. In particular, ${\rm Spec}(L_G \otimes \Q) \cong Y_G$. 
        \item If $c > 1$, the rings $L_G$ and $\Z[Y_G]$ are not isomorphic.
    \end{enumerate} 
\end{lemma}
\begin{proof}
Set $R=\Z[Y_G]$. We first prove (1). Since
\[
    \Phi(u')\Phi(v')=(yz-A^*(x))(xz-B^*(y))=cz+A^*(x)B^*(y),
\]
the element $z$ belongs to the image of $\Phi[1/c]$. Hence $\Phi[1/c]$ is surjective. Its inverse sends
\[
    x\longmapsto u,\qquad y\longmapsto v,\qquad
    z\longmapsto \frac{u'v'-A^*(u)B^*(v)}{c}.
\]
This is well-defined since
\[
    uv\bigl(u'v'-A^*(u)B^*(v)\bigr)
    =A(u)B(v)-(A(u)-c)(B(v)-c)=c(A(u)+B(v)-c),
\]
is the defining equation of $R[1/c]$. This proves (1).

Suppose now that $c > 1$, and choose a prime $p$ dividing $c$.
Set $k=\F_p$, write $a,b$ for the reductions of
$A,B$ modulo $p$, and consider
\[
    R_p:=R\otimes_\Z k
       \cong k[x,y,z]/(xyz-a(x)-b(y))
\]
and 
\[
    (L_G)_p:=L_G\otimes_\Z k
    \cong k[u,u',v,v']/(uu'-b(v),\, vv'-a(u)).
\]
Set $J = (uu' - b(v), vv' - a(u))$, and notice that $J$ has height $2$ in $k[u,u',v,v']$. Finally, set $t = \gcd(xy,a(x)+b(y))$. 
We consider four cases. 

{\bf Case 1.} $t = 1$.
Then $a,b$ are both nonzero, and $R_p$ has a unique minimal prime. On the other hand, we claim that $(L_G)_p$ has at least two. The ideal $(u,v)$ is prime, has height 2 in $k[u,u',v,v']$ and contains $J$; it is thus a minimal prime of $(L_G)_p$. Consider the map $\phi: k[u,u',v,v'] \to k[u^{\pm 1}, v^{\pm 1}]$ given by fixing $u,v$ and where $u' \mapsto \frac{b(v)}{u}$ and $v ' \mapsto \frac{a(u)}{v}$. Then $T:= \phi^{-1}((0))$ is prime, of height 2 and contains $J$, so is also a minimal prime of $(L_G)_p$. To see they are distinct, note $uv \in (u,v)$ but $uv \not\in T$. 

{\bf Case 2.} $t = x$. Then  $b = 0$ and $a \neq 0$. Then $xyz -a(x) = x(yz - a^*(x))$, and so $R_p$ has two minimal primes; here $a^*$ denotes the reduction of $A^*$ modulo $p$. On the other hand, we claim $(L_G)_p$ has at least three. It is clear that $(u,v)$ and $(u,v')$ are distinct minimal primes. Consider the map $\phi: k[u,u',v,v'] \to k[u^{\pm 1}, v^{\pm 1}]$ given by fixing $u,v$, $u' \mapsto 0$ and $v' \mapsto \frac{a(u)}{v}$. Then $\phi^{-1}(0)$ is prime, of height 2, contains $J$ and hence is minimal in $(L_G)_p$. It is distinct from $(u,v)$ and $(u,v')$ since $u \not\in \phi^{-1}(0)$.  

{\bf Case 3.} $t = y$. Then $a = 0$ and $b \neq 0$. This case is similar to Case 2, and we omit it.

{\bf Case 4.} $t = xy$. Then $a = b = 0$. Then $R_p$ has three minimal primes, and $(L_G)_p$ has four. 
 
In all cases, $R_p$ and $(L_G)_p$ have distinct numbers of minimal primes. Therefore, $R$ and $L_G$ cannot be isomorphic.
\end{proof}

\section{Isomorphisms of Mordell--Schinzel surfaces}\label{sec:companion}
In this section, we construct
isomorphisms between
pairs of Mordell--Schinzel surfaces
by rescaling coordinates and using the exchange relations.
We also give conditions
under which these isomorphisms
preserve positive integral points.

\subsection{Coordinate rescaling}\label{subsec:coord-rescale}
Over a field $K$ of characteristic zero, the substitution $x_i = \lambda_ix_i'\; (0 \leq i \leq d)$, with $\lambda_i \in K^\times$, identifies $Y_H$ with $Y_G$, where   
\begin{equation}\label{eq:coordinate-rescaling}
     H(x_1', \ldots, x_d') = (\lambda_0\cdots \lambda_d)^{-1} G(\lambda_1x_1',\ldots, \lambda_dx_d').
\end{equation}
   
\begin{example}\label{ex:1-1-single-isom-class}
    For $G = ax + by + c$ with $abc \neq 0$, \eqref{eq:coordinate-rescaling} gives 
    \[
        Y_G \cong Y_{x+y+1}, \quad (x,y,z) \mapsto (ax/c, by/c, cz/(ab)).
    \] 
\end{example}

\begin{remark}\label{rem:Q-MS-surface-has-model}
    If $G$ has rational coefficients, and $t$ denotes the least common multiple of the denominators of all coefficients of $G$, then $(x_0,\ldots, x_d) \mapsto (tx_0,\ldots, x_d)$ identifies $Y_G$ with the integral model $Y_{tG}$. In particular, every Mordell--Schinzel variety defined over $\Q$ has an integral model.
\end{remark} 
\subsection{Two isomorphisms}\label{subsec:two-isoms}
In this subsection, we fix a pure $G\in K[x,y]$ with $m,n\geq1, c\neq0$ and decomposition \eqref{eq:pure-poly}.
Recall that $a$ and $b$
denote the leading coefficients of $A$ and $B$ respectively.

The following construction is due to
Koll\'{a}r--Villalobos-Paz~\cite[Paragraph~26]{KVP:cubic},
and has also appeared in
\cite[Equations~(52)--(53)]{schinzel-2015}. Set
\[
    r_1 = ac^{m - 2}, \quad r_2 = bc^{n - 2},
\]
and define
\begin{equation}\label{eq:reciprocal-tilde}
    A^\dagger(x) = \frac{x^m}{c}A(c/x), \qquad
    B^\dagger(y) = \frac{y^n}{c}B(c/y).
\end{equation}
Then set
\begin{equation}\label{eq:G_1}
    \begin{aligned}
        G_1 &= A^\dagger(x) + r_1B(y) - r_1c, \\
        G_2 &= r_2A(x) + B^\dagger(y) - r_2c.
    \end{aligned}
\end{equation}
By construction, $A^\dagger, B^\dagger \in K[t]$ are monic, of constant term $r_1c$ and $r_2c$ respectively. Their roots, counted with multiplicities, are $c/\rho$ as $\rho$ runs through the roots, also counted with multiplicities, of $A$
and $B$ respectively. Recall from \eqref{eq:A-B-*} that $A^*$ and $B^*$ are defined by $A(x) - c = xA^*(x)$ and $B(y) - c = yB^*(y)$.

\begin{proposition}[{\cite[\S 26]{KVP:cubic}}]
\label{prop:birational-are-isom}
    Set $u=yz-A^*(x)$ and $v=xz-B^*(y)$. The maps
    \begin{align}\label{eq:mu_1}
        &\mu_{G,1}(x,y,z) = \left(u,y,\frac{(A^\dagger(u)-r_1c)/u+r_1x}{y}\right), \\\label{eq:mu_2}
        &\mu_{G,2}(x,y,z) = \left(x,v,\frac{(B^\dagger(v)-r_2c)/v + r_2y}{x} \right), 
    \end{align}
    define $K$-isomorphisms $Y_G \to Y_{G_1}$ and $Y_G \to Y_{G_2}$ respectively. On $D(xy)$, their first two coordinates are
    $(B(y)/x, y)$ and $(x, A(x)/y)$ respectively.

\end{proposition}
\begin{proof}
It suffices to treat $\mu_{G,1}$. We first show that the image of $\mu_{G,1}$ lies in $Y_{G_1}$. Set $\mu_{G,1}(x,y,z) = (u,y,w)$ and
$Q(t)=(A^\dagger(t)-r_1c)/t$. Then 
\[
    uyw = A^\dagger(u) +r_1ux - r_1c = A^\dagger(u) + r_1B(y)-r_1c = G_1(u,y).
\]
On $D(xy)$, $z = (A(x)+B(y)-c)/(xy)$ and so $\mu_{G,1}$ is determined by its first two coordinates. The last claim follows from the identity $yz - A^*(x) = B(y)/x$. It remains to prove that $\mu_{G,1}$ is regular, with regular inverse. By construction, $A^\dagger(u)$ has constant term $r_1c$ and degree at least 1, so $Q(u)$ is regular. It remains to show that $Q(u) + r_1x$ is divisible by $y$ in $K[Y_G]$. Note that $xu=B(y)$, and thus $xu\equiv c \mod y$ and $A(x) \equiv 0 \mod y$. Since $c \neq 0$, $u$ is a unit modulo $y$, and so 
\[
 A^\dagger(u)=\frac{u^m}{c}A(c/u) \equiv \frac{u^m}{c}A(x) \equiv 0 \mod y,
\]
and thus $u(Q(u)+r_1x) \equiv A^\dagger(u) -r_1c + r_1xu \equiv 0 \mod y $, and since $u$ is a unit we have that $Q(u)+r_1x \equiv 0 \mod y$, i.e. $w \in K[Y_G]$. Its inverse is
\[
 x=\frac{yw-Q(u)}{r_1},\qquad y=y,\qquad
 z=\frac{A^*(x)+u}{y}.
\]
Indeed, on $Y_{G_1}$ we again have $xu=B(y)$. Modulo $y$,
$xu=c$ and $A^\dagger(u)=0$, so $A(x)=0$ and $A^*(x)+u=0$.
The inverse is therefore regular. The two compositions are the
identity on $D(xy)$, hence everywhere. Interchanging $A,x$ with
$B,y$ proves the second assertion.
\end{proof}

\begin{remark}
    The morphism $\mu_{G,i}$ of Proposition \ref{prop:birational-are-isom} agrees with that in \cite[26]{KVP:cubic}, up to a coefficient scaling by a factor of $r_i$. More precisely, our definition of $G_1$ and $G_2$ ensures $G_1(x,0) = A^\dagger(x)$ and $G_2(0,y) = B^\dagger(y)$ are monic and with constant terms $r_1c$ and $r_2c$, whereas \cite[26]{KVP:cubic} require $G_1$ and $G_2$ have same constant term as $G$. 
\end{remark}

\begin{example}\label{ex:2-2-isoms}
    Fix $g_0 = ax^2 + by^2 + c$ with $a,b,c \in K^\times$. Set $M = ab$, and define $g_k \in K[x,y]$ for each $k \geq 1$ by
    \[
    g_k = \begin{cases}
        x^2 + My^2 + aM^{k-1}c, &\text{if } k \text{ odd}, \\
        Mx^2 + y^2 + aM^{k-1}c&\text{if } k \text{ even,}
    \end{cases}
    \]
    Then, 
    \[
        Y_{g_0} \stackrel{\mu_{g_0,1}}{\longrightarrow} Y_{g_1} \stackrel{\mu_{g_1,2}}{\longrightarrow} Y_{g_2} \stackrel{\mu_{g_2,1}}{\longrightarrow} Y_{g_3}
        \stackrel{\mu_{g_3,2}}{\longrightarrow} Y_{g_4} \stackrel{\mu_{g_4,1}}{\longrightarrow}\cdots
    \]
    where $\mu_{g_0,1}(x,y,z) =(yz - ax,y,z)$, and for $k \geq 1$, 
    \[
        \mu_{g_k,i}(x,y,z) = \begin{cases}
        (yz-abx,y,z)&\text{if } i = 1 \text{ and } k \text{ even,} \\
        (x,xz-aby,z), &\text{if } i = 2 \text{ and } k \text{ odd.}
    \end{cases}
    \] 
\end{example}
 
 The maps $\mu_{G,1}$ and $\mu_{G,2}$ define automorphisms of $Y_G$ precisely when $r_1=r_2 = 1$, $A^\dagger = A$ and $B^\dagger = B$. 

 \begin{example}\label{ex:2-2-stable}
     Suppose $G = x^2 + a_1x+ y^2 + b_1y + c$ with $c \neq 0$. Then $r_1 = r_2 = 1$, $A^\dagger=A$ and $B^\dagger = B$ so $\mu_{G,1}$ and $\mu_{G,2}$ are automorphisms of $Y_G$. The case $c= 1$ was studied by Mills \cite{mills}, and the case $a_1 = b_1 = 0$ by Barnes \cite{barnes}.
 \end{example}

\subsection{Companion surfaces}\label{subsec:companion}
In this subsection, we fix a pure $G \in \Znn[x, y]$ such that
$mn \geq 4$, $c \neq 0$ and $abc^{m+n-4}=1$. The last assumption is equivalent to the assumption that one of the following holds:
\begin{align}
    &a = b = c = 1;
    \label{eq:am-bn-c=1}\\
    &m = n = 2 \text{ and } a = b = 1.
    \label{eq:m-n=2-a-b=1}
\end{align}
In either case, $r_1 = r_2 = 1$ and $A^\dagger(0) = B^\dagger(0) = c$,
so substituting into \eqref{eq:G_1} gives
\begin{align*}
    G_1 &= A^\dagger(x) + B(y) -c, \qquad A^\dagger(0) = B(0) = c, \\
    G_2 &= A(x) + B^\dagger(y) -c,  \qquad  A(0) = B^\dagger(0) = c.
\end{align*}
Note that $(A^\dagger)^\dagger = A$ and $(B^\dagger)^\dagger = B$. Thus, repeating the above for $G_1$ and $G_2$ instead of $G$ yields a single new polynomial
\[
    G_3 = A^\dagger(x) + B^\dagger(y) - c, \qquad A^\dagger(0) = B^\dagger(0) = c.
\]
For $R \in \{A,A^\dagger\}$, $S \in \{B,B^\dagger\}$, we have $R(0) = S(0) =c$, thus the decomposition of $G_{R,S} = R(x) + S(y) -c$ is of the form \eqref{eq:pure-poly};
we write $Y_{R,S} := Y_{G_{R,S}}$ and $\mu_{R,S,i} := \mu_{G_{R,S},i}$ for
\[
    \begin{tikzcd}[row sep = tiny]
        \mu_{G_{R,S},1}: Y_{R,S} \arrow[r] & Y_{R^\dagger,S} \\
        \mu_{G_{R,S},2}: Y_{R,S} \arrow[r] & Y_{R,S^\dagger}.
    \end{tikzcd}
\]
We call $Y_{A,B}$, $Y_{A^\dagger,B}$, $Y_{A,B^\dagger}$ and
$Y_{A^\dagger,B^\dagger}$ the {\it companion surfaces} of
$Y_{A,B}$, identifying any that coincide. These are the objects
of the groupoid $\mathcal{G}_{A,B}$ defined in \S \ref{subsec:outline};
its morphisms are the isomorphisms generated by
the maps $\mu_{R,S,i}$ and their inverses.

\begin{remark}
    If $G$ satisfies \eqref{eq:am-bn-c=1} then $c = 1$ and so $A^\dagger(x) = x^m A(1/x)$, $B^\dagger(y) = y^n B(1/y)$. If $G$ satisfies \eqref{eq:m-n=2-a-b=1}, then $A^\dagger = A$ and $B^\dagger= B$.
\end{remark}

\begin{proposition}
    \label{prop:integral-normalization}
    For $R \in \{A,A^\dagger\}$, $S \in \{B,B^\dagger\}$,
    the maps $\mu_{R,S,1}$ and $\mu_{R,S,2}$ are isomorphisms over $\Z$
    with inverses $\mu_{R^\dagger,S,1}$ and $\mu_{R,S^\dagger,2}$, respectively.
    Furthermore, they induce bijections on positive integral points. 
\end{proposition}

\begin{proof}
    If $a=b=c=1$, the proof of
    Proposition~\ref{prop:birational-are-isom} applies over $\Z$:
    $r_1=r_2=c=1$, and all the displayed polynomials have integral
    coefficients. Since $(R^\dagger)^\dagger=R$ and
    $(S^\dagger)^\dagger=S$, the stated maps are inverse.
    If $m=n=2$ and $a=b=1$, write
    $R(t)=t^2+a_1t+c$ and $S(t)=t^2+b_1t+c$.
    Then the maps $\mu_{R,S,i}$ send
    \[
        \begin{tikzcd}[row sep = tiny]
            \mu_{R,S,1}: (x,y,z) \arrow[r, mapsto] & \paren{yz-x-a_1,y,z} \\
            \mu_{R,S,2}: (x,y,z) \arrow[r, mapsto] & \paren{x,xz-y-b_1,z}
        \end{tikzcd}
    \] 
    and both maps are involutions over $\Z$.
    The positivity statement follows from the last assertion in Proposition~\ref{prop:birational-are-isom};
    the exchanged coordinates $S(y)/x$ and $R(x)/y$ are positive
    and so the defining equation of each target surface then makes
    the third coordinate positive as well.
    The same argument applies to the inverse maps
    $\mu_{R^\dagger,S,i}$.
\end{proof}
Proposition \ref{prop:integral-normalization} may fail if we do not assume \eqref{eq:am-bn-c=1} or \eqref{eq:m-n=2-a-b=1} holds, as the following example shows.

\begin{example}
    \label{eq:2-2-rational}
    In the notation of Example \ref{ex:2-2-isoms}, the two maps
    \[
        \begin{tikzcd}[row sep = tiny]
            \mu_{g_0,1}^{-1}: (x,y,z) \arrow[r, mapsto] & \paren{\frac{yz-x}{a},y,z} \\
            \mu_{g_1,2}^{-1}: (x,y,z) \arrow[r, mapsto] & \paren{x, \frac{xz-y}{ab}, z}
        \end{tikzcd}
    \] 
    are defined over $\Z$ only if $a,b \in \{\pm 1\}$. 
\end{example}

The {\it vertex group} of $Y_{A,B}$, denoted $\Gamma_{A,B}$, is the group of morphisms of $\mathcal{G}_{A,B}$ that are automorphisms of $Y_{A,B}$. When $G \in \Znn[x,y]$, it follows from Proposition \ref{prop:integral-normalization} that $\Gamma_{A,B}$ is a subgroup of ${\rm Aut}_\Z^+(Y_{A,B})$.
\begin{example}
    If $A^\dagger = A$ and $B^\dagger = B$, then 
    \[
        \Gamma_{A,B} = \langle \mu_{A,B,1}, \mu_{A,B,2}\rangle.
    \]
    In particular, $\Gamma_{A,B} \cong D_\infty = \langle \sigma, \mu : \sigma^2 = \mu^2 = 1 \rangle$ whenever $mn \geq 4$. This is the case for instance when $G$ is as in Example \ref{ex:2-2-stable}.
\end{example}

When $m = n = 3$, condition \eqref{eq:am-bn-c=1} implies $a = b = c = 1$. It was shown in \cite[Theorem 28]{KVP:cubic} that $\Gamma_{A,B}$ has finite index in ${\rm Aut}_\Z(Y_G)$. 

When $A^\dagger = A$, $B^\dagger = B$ and $A$ and $B$ have only simple roots (say, over $\C$), the surface $Y_{A,B}$ is smooth (cf. Lemma \ref{lem:singular-locus-general}). In this case, the techniques of \cite{blanc-dolgachev} together with Lemma \ref{lem:2-integral-forms} can be used to show that $\Gamma_{A,B}$ has finite index in, and is often equal to, ${\rm Aut}_\Z(Y_{A,B})$. Since we will not be needing this result, we omit its proof.

\section{Generalised upper cluster algebras and families of Mordell--Schinzel surfaces}\label{sec:gca}
Fix $m, n \geq 1$.
In this section, we identify the coordinate ring of
the coefficient family $\calm_{m,n}$
with a generalised upper cluster algebra $\calu_{m, n}$,
and then establish a correspondence between
adapted friezes of $\calu_{m, n}$
and positive integral points on Mordell--Schinzel surfaces.

\subsection{Coefficient families and seed patterns}\label{sec:ca}
Let ${\bf s}=(a_1,\ldots,a_m,b_1,\ldots,b_n,c)$ be a tuple of formal variables, set 
\[
A(t) = c + \sum_{i = 1}^m a_it^i, \quad B(t) = c + \sum_{j = 1}^n b_jt^j, \quad  \LL = \Z[{\bf s}], \quad \LL_+ = \Znn[{\bf s}],
\]
and consider the flat family of surfaces $
\pi: \calm_{m,n} \longrightarrow \calb_{m,n}$
where $\calm_{m,n}$ is the
hypersurface given by $xyz = A(x) + B(y) - c$ in $\A^3_\Z \times \calb_{m,n}$, $\calb_{m,n} = \operatorname{Spec}(\LL)$ and $\pi$ is the restriction of the canonical projection onto $\calb_{m,n}$. Finally, set
\[
    \calb_{m,n}^+ =
    \set{(a_1, \ldots, a_m, b_1, \ldots, b_n, c)
    \in \Znn^{m + n + 1} : a_mb_nc > 0}.
\]
Recall the notation of Appendix \ref{sec:appendix}. Let $ \P = {\rm Trop}({\bf s})$ and $\F = \Q\P(x_0,x_1)$, with $x_0$ and $x_1$ algebraically independent over $\Q\P$. Denote by $[\Sigma_{m,n}]$ the unique seed pattern of geometric type defined by the initial seed   
\[
    \Sigma_{m,n} =  ((x_0, x_1), ({\bf p}_1, {\bf p}_2), \begin{psmallmatrix}0 & m \\ -n & 0\end{psmallmatrix}),
\]
where ${\bf p}_1 = (b_n, b_{n - 1}, \ldots, b_1,c)$ and  ${\bf p}_2 = (c, a_1,a_2, \ldots, a_m)$. Denote by $\{x_k : k \in \Z\}$ and by $\{E_k : k \in \Z\}$ the set of cluster variables and exchange polynomials of $[\Sigma_{m,n}]$ respectively. Thus
\begin{equation}\label{eq:E0-E1}
 E_0(t)=A(t),\qquad E_1(t)=B(t),
\end{equation}
and
\begin{equation}\label{eq:cv-recursive}
 x_kx_{k+2}=E_{k+1}(x_{k+1})\qquad(k\in\Z).
\end{equation}
Put $d_k=m$ for even $k$ and $d_k=n$ for odd $k$, and write $E_k(t)=\sum_{r=0}^{d_k}p_{k;r}t^r$. Then 
\begin{align}
 E_{k+2}(t)&=\frac{t^{d_k}E_k(E_{k+1}(0)/t)}
 {\bigoplus_{\ell=0}^{d_k}p_{k;\ell}E_{k+1}(0)^\ell}
 &&(k\geq0),\label{eq:Ek+}\\
 E_k(t)&=\frac{t^{d_k}E_{k+2}(E_{k+1}(0)/t)}
 {\bigoplus_{\ell=0}^{d_k}p_{k+2;\ell}E_{k+1}(0)^\ell}
 &&(k<0).\label{eq:Ek-}
\end{align}
In particular, $E_k\in\LL_+[t]$. Set 
\[
 \calu_{m,n}=\bigcap_{k\in\Z}L_k,
 \qquad \calu_{m,n}^+=\bigcap_{k\in\Z}L_k^+,
\]
where the intersection is taken in $\F$ and 
\begin{equation}\label{eq:Laurent-chart}
 L_k=\LL[x_k^{\pm1},x_{k+1}^{\pm1}],\qquad
 L_k^+=\LL_+[x_k^{\pm1},x_{k+1}^{\pm1}].
\end{equation}
Thus $\calu_{m,n}$ is a generalised upper cluster algebra
of geometric type, rank $2$,
and with non-inverted frozen variables.
The strong positive Laurent phenomenon (Theorem \ref{thm:appendix-strong-pos-laur-phen}) gives $x_j \in \calu_{m,n}^+$ for all $j \in \Z$. 

\begin{example}
    When $m = n = 2$, set
    \[
        E_0(t) = at^2 + a_1t + c, \quad E_1(t) = bt^2 + b_1t + c.
    \]
    Then, 
    \begin{equation*}
        E_2(t)=t^2+a_1t+ac, \quad E_3(t) =t^2+ab_1t+aMc, \quad E_4(t)=t^2+a_1Mt+aM^2c,
    \end{equation*}
    where $M = ab$ (compare with Example \ref{ex:2-2-isoms}). Set $x=x_0$ and $y=x_1$. Substitution in \eqref{eq:cv-recursive} gives
\begin{align*}
 x_2&=\frac{by^2+b_1y+c}{x}, \quad x_3=\frac{B(y)^2+a_1xB(y)+acx^2}{x^2y},\\
 x_4&=\frac{1}{x^3y^2}\Bigl(B(y)^3+2a_1xB(y)^2 +(a_1^2+ab_1y+2ac)x^2B(y)\\
    &\qquad \qquad +a_1a(b_1y+2c)x^3+a^2cx^4\Bigr).
\end{align*}

\end{example}

\begin{example}\label{ex:G2}
For the seed pattern $[\Sigma_{3, 1}]$, write $b = b_1$.
We have
    $x_{k + 8} = x_k$ and $E_{k + 8} = E_k$ for all $k \in \Z$.
    The polynomials $E_0, \ldots, E_7$ are given in
    Figure~\ref{fig:G2-exchange}.
    \begin{figure}[H]
        \centering
        \begingroup
            \renewcommand{\arraystretch}{1.25}
            \[
                \begin{array}{c@{\qquad}l}
                    \hline
                    k & \displaystyle E_k(t) \\
                    \hline
                    0 & \displaystyle
                    c + a_1t + a_2t^2 + a_3t^3 \\
                    1 & \displaystyle c + bt \\
                    2 & \displaystyle a_3c^2 +
                    a_2c
                    t + a_1t^2 + t^3 \\
                    3 & \displaystyle a_3bc + t \\
                    4 & \displaystyle
                    a_3^2b^3c +
                    a_1a_3b^2t + 
                    a_2bt^2 + t^3
                    \\
                    5 & \displaystyle a_3b^2 + t \\
                    6 & \displaystyle
                    a_3b^3 + a_2b^2t +
                    a_1bt^2 +
                    c
                    t^3 \\
                    7 & \displaystyle b + t \\
                    \hline
                \end{array}
            \]
        \endgroup
        \caption{The exchange polynomials of $[\Sigma_{3, 1}]$}
        \label{fig:G2-exchange}
    \end{figure}

    Set $x=x_0$, $y=x_1$, and
    \begin{align*}
        \Delta_1 & = a_3c^2x^3 + a_2bc x^2y
        + a_1b^2xy^2 + b^3y^3 + a_2c^2x^2
        \\
        & \quad + 2a_1bc xy + 3b^2c y^2
        + a_1c^2x + 3bc^2y + c^3,
        \\
        \Delta_2 & = a_3^2c x^6 + a_2^2c x^4 +
        b^3y^3
        +3b^2c y^2 + 3bc^2y + c^3
        \\
        & \quad + (bx^2y + c x^2)a_1^2
        +2(b^2xy^2 + 2bc xy + c^2x)a_1
        \\
        & \quad + \bigl(2a_2c x^5 + 3bc
        x^3y + 2c^2x^3
        +(bx^4y + 2c x^4)a_1\bigr)a_3
        \\
        & \quad + \bigl(b^2x^2y^2 + 3bc
        x^2y + 2c^2x^2
        +(bx^3y + 2c x^3)a_1\bigr)a_2,
        \\
        \Delta_3 & = a_3c x^3 + a_2c
        x^2 + b^2y^2 + 2bc y
        \\
        & \quad + (bxy + c x)a_1 + c^2.
    \end{align*}
    As rational functions of $x$ and $y$, the cluster variables are

        \begin{align}
            x_2 &= \frac{by + c}{x},  \quad            x_7 = \frac{a_3x^3 + a_2x^2 + a_1x +
            c}{y},
            \\[1.2ex]
            x_6 &= \frac{a_3x^3 + a_2x^2 + a_1x +
            by +
            c}{xy},\label{eq:z-in-g2}
        \end{align}
    and $x_3 = \frac{\Delta_1}{x^3y}, x_4 = \frac{\Delta_3}{x^2y}, 
            x_5 = \frac{\Delta_2}{x^3y^2}$.

\end{example}

\begin{theorem}
    \label{thm:integral-form-uca}
    The $\LL$-algebra homomorphism $\Psi: \Z[\calm_{m, n}] \longrightarrow \calu_{m,n}$ given by 
\[
            (x,y,z)\mapsto \paren{x_0, x_1, \dfrac{A(x_0) + B(x_1) - c}{x_0x_1}}
    \]
    is an isomorphism, and $\Psi(z) \in \calu_{m, n}^+$.
\end{theorem}
\begin{proof}
The map $\Psi: \Z[\calm_{m,n}] \to \F$ is well-defined, since the ideal $(xyz - A(x)-B(y) + c)$ is contained in $\ker\Psi$. To show its codomain is $\calu_{m,n}$ it suffices to prove that $\Psi(z) \in \calu^+_{m,n}$. 

Let $[\widetilde{\Sigma}_{m,n}] = \{\Sigma_k : k \in \Z\}$ be a (non-normalised) seed pattern of rank $2$ such that $\Sigma_0 = \widetilde \Sigma_{m,n} = ((x_0,x_1), (\widetilde{\bf p}_1,\widetilde{\bf p}_2),M )$ where
\[
    \widetilde p_{i,r} = \frac{p_{i,r}}{c}, \quad M = \begin{pmatrix}
        0 & m \\ -n & 0
    \end{pmatrix}, \qquad i \in \{1,2\}, \; 0 \leq r \leq d_i.
\]
We emphasise that $\widetilde \Sigma_{m,n}$ is not normalised, and thus $[\widetilde{\Sigma}_{m,n}]$ is not uniquely determined by $\widetilde \Sigma_{m,n}$. Denote by $\{\widetilde x_k : k \in \Z\}$ the set of cluster variables of $[\widetilde \Sigma_{m,n}]$, with $\widetilde x_0 = x_0$ and $\widetilde x_1 = x_1$. It follows by an inductive argument on $\T_2$ (see also \cite[Proposition~7.13 and Remark~7.15]{burcroff-Lee-mou:generalised}) that 
\[
    \frac{\widetilde x_k}{x_k} \in \P, \qquad \forall k \in \Z,
\]
and so $\calu^+([\widetilde\Sigma_{m,n}]) = \calu^+([\Sigma_{m,n}])$, where
\[
\calu^+([\widetilde\Sigma_{m,n}]) = \bigcap_{k\in \Z} \Znn\P [\widetilde x_k^{\pm 1}, \widetilde x_{k+1}^{\pm 1}] \quad \text{and} \quad \calu^+([\Sigma_{m,n}]) = \bigcap_{k\in \Z} \Znn\P [ x_k^{\pm 1}, x_{k+1}^{\pm 1}].
\]
Since $\widetilde p_{1,n} = \widetilde p_{2,0} = 1$, Example \ref{ex:non-normalised-greedy} shows that the greedy element 
\[
    x[1,1]
    = \frac{1+\sum_{i=1}^m\widetilde p_{2,i}(\widetilde x_0)^i +\sum_{j=1}^n\widetilde p_{1,n-j}(\widetilde x_1)^j}
        {\widetilde x_0 \widetilde x_1} = \frac{1+\sum_{i=1}^m \frac{a_i}{c} x_0^i +\sum_{j=1}^n \frac{b_j}{c} x_1^j}
        { x_0x_1}
\]
belongs to $\calu^+([\Sigma_{m,n}])$. Thus $\Psi(z) = cx[1,1]$ also belongs to $\calu^+([\Sigma_{m,n}])$. It remains to show that $\Psi(z) \in \LL_+[x_k^{\pm 1}, x_{k+1}^{\pm 1}] \subset \Znn\P[x_k^{\pm 1}, x_{k+1}^{\pm 1}]$ for each $k \in \Z$. To see this, we apply the semifield homomorphism of Lemma~\ref{lem:sf-homo} (which maps every $x_j$ to $1$ and fixes the generators of $\LL_+$) to $\Psi(z)$: 
\begin{align*}
 \Psi(z)\big\vert_{\P}&=\frac{\sum_{i=1}^m a_ix_0^i + \sum_{j=1}^n b_jx_1^j + c}{x_0x_1} \big\vert_{\P} \\
 &= a_1\oplus \cdots\oplus a_m
                \oplus b_1\oplus\cdots\oplus b_n\oplus c=1.
\end{align*}
Therefore, $\Psi(z) \in \calu_{m,n}^+$. 

We now prove that $\Psi$ is an isomorphism. It is clear that $\Psi$ is injective. It remains to show that $\Psi$ is surjective. 
Identify $\Z[\calm_{m,n}]$ with its image and write $x=x_0$, $y=x_1$, $z=\Psi(z)$. Also write $x' = x_2$ and $y' = x_{-1}$. The relations \eqref{eq:cv-recursive} are
\begin{equation}\label{eq:one-step-reduction}
    y'=xz-B^*(y),\qquad x'=yz-A^*(x).
\end{equation}
Thus
\[
 L_{-1}=\Z[\calm_{m,n}][(xy')^{-1}],\quad
 L_0=\Z[\calm_{m,n}][(xy)^{-1}],\quad
 L_1=\Z[\calm_{m,n}][(yx')^{-1}].
\]
Since \(\mathcal U_{m,n}\subseteq L_{-1}\cap L_0\cap L_1\), it is enough to
show
\[
    L_{-1}\cap L_0\cap L_1=\Z[\calm_{m,n}].
\]
The ideals $(x,y)$, $(x,x')$ and $(y,y')$ have height two. Indeed, the relations \eqref{eq:one-step-reduction} are linear in $a_1$ and in $b_1$, respectively. No height-one prime of $\Z[\calm_{m,n}]$
can therefore contain all three products $xy'$, $xy$, $yx'$.
Since $\Z[\calm_{m,n}]$ is normal, and a normal noetherian domain is the intersection of its
localizations at height-one primes,
\[
 \calu_{m,n}\subseteq L_{-1}\cap L_0\cap L_1=\Z[\calm_{m,n}].
\]
This concludes the proof of the theorem.
\end{proof}

\subsection{Friezes}
\label{sec:friezes}
Friezes and their combinatorial counterparts, frieze patterns, have a long history of study, see e.g. \cite{coxeter,coxeter-conway-2,caldero-chapoton,assem-reutenauer-smith,morier-genoud}. The following definition of a frieze is essentially \cite[Definition 3.1]{dhl}.
\begin{definition}
    \label{def:frieze}
    For $p \in \calb_{m,n}^+$, a frieze adapted to $p$ is a ring homomorphism $f : \calu_{m, n} \longrightarrow \Z$ whose restriction to $\LL$ corresponds to evaluation at $p$ and for which $f(x_i) >0$ for every $i \in \Z$. Write ${\rm Frieze}(\calu_{m,n}, p)$ for the set of these homomorphisms.
\end{definition}

Fix a pure $G \in \Znn[x, y]$ of bidegree $(m, n)$ with
$c > 0$, decomposed uniquely as $G = A(x) + B(y) - c$. Write $A,B$ as in Section \ref{sec:ca}, and set
$p_G = (a_1, \ldots, a_m, b_1, \ldots, b_n, c)
\in \calb_{m,n}^+$.

\begin{theorem}\label{thm:frieze-and-integer-points}
    The evaluation map
    \[
    \Frieze(\calu_{m, n}, p_G) \longrightarrow Y_G(\Zpos), \quad f \mapsto \bigl(f(x_0), f(x_1), f(\Psi(z))\bigr)
    \]
    is a bijection.
\end{theorem}
\begin{proof}
By Theorem~\ref{thm:integral-form-uca}, a point of
$Y_G(\Zpos)$ determines a unique homomorphism
$f:\calu_{m,n}\to\Z$ extending evaluation at $p_G$ and taking
$x_0,x_1,\Psi(z)$ to its coordinates. The proof of that theorem
shows that every $x_j$ has a Laurent expansion in $x_0,x_1$ with
coefficients in $\LL_+$, so $f(x_j)\geq0$.
The constant and leading coefficients of each $E_j$ are monomials
in $a_m,b_n,c$, as follows inductively from
\eqref{eq:Ek+}--\eqref{eq:Ek-}. In particular, its constant term
specialises to a positive integer. Thus
\[
 f(x_{j-1})f(x_{j+1})=f(E_j(x_j))>0
 \qquad(j\in\Z),
\]
and every $f(x_j)$ is positive. Hence $f$ is an adapted frieze.
Conversely, an adapted frieze has $f(x_0),f(x_1)>0$ and
\[
 f(\Psi(z))=
 \frac{G(f(x_0),f(x_1))}{f(x_0)f(x_1)}>0.
\]
The defining equation gives a point of $Y_G(\Zpos)$.
The constructions are inverse because $x_0,x_1,\Psi(z)$ generate
$\calu_{m,n}$ over $\LL$.
\end{proof}

\section{Finiteness and density of positive integral points}
\label{sec:theoremA}
In this section, we classify the finiteness of $Y_G(\Zpos)$
for both $c > 0$ and $c = 0$,
and then prove Zariski density when
$c > 0$ and $mn > 4$.

For Sections~\ref{subsec:finite} and~\ref{ss:infinite-type},
fix a pure $G \in \Znn[x, y]$ with $m, n \geq 0$
and $c > 0$.

\subsection{Positive constant term: finite type}\label{subsec:finite}

\begin{proposition}\label{prop:finite-case}
    Suppose that $mn < 4$. Then $Y_G(\Zpos)$ is finite.
\end{proposition}

If $m = 0$, then $xyz = B(y)$, so $y \mid c$. For each
such positive $y$, the positive integers $x, z$ satisfy
$xz = B(y)/y$, giving only finitely many possibilities.
The case $n = 0$ is symmetric. Hence, assume $m, n \geq 1$.
Up to swapping $x$ and $y$, we have
\begin{equation}\label{eq:Gp-finite-type}
    G = a_3x^3 + a_2x^2 + a_1x + by + c,
\end{equation}
where $a_1, a_2, a_3 \in \Znn$, not all zero, and $b, c \in \Zpos$. 

The most direct approach to proving Proposition \ref{prop:finite-case} is as follows. Set $i = \max(i : a_i \neq 0)$ and write $p_i = (a_1,\ldots, a_i,b,c) \in \calb_{i,1}^+$. By Theorem \ref{thm:frieze-and-integer-points}, $Y_G(\Zpos)$ is in bijection with ${\rm Frieze}(\calu_{i,1}, p_i)$. Thus one is reduced to showing that ${\rm Frieze}(\calu_{i,1}, p_i)$ is a finite set, in each of the cases $i = 1,2,3$. This is somewhat tedious. 

We therefore take a more uniform approach. Define 
\[
    p_G = (a_1,a_2,a_3,b,c) \in \calb_{3,1}.
\]
We denote by $\{x_i : i \in \{0,\ldots, 7\}\}$ the set of cluster variables of $\calu_{3,1}$ computed in Example \ref{ex:G2}, and by ${\rm ev}_{p_G}$ the restriction to $\calu_{3,1}$ of the evaluation map 
\[
    \Z[\calb_{3,1}][x_0^{\pm 1},x_1^{\pm 1}] \longrightarrow \Z[x^{\pm 1}, y^{\pm 1}], \qquad x_0 \mapsto x, \; x_1 \mapsto y,
\]
which specialises the coordinates of $\calb_{3,1}$ to their numerical values prescribed by $p_G$. For each $i = 0,\ldots, 7$, set $u_i = {\rm ev}_{p_G} x_i$. Note that each $u_i= u_i(x,y)$ is a Laurent polynomial in $x = x_0, y = x_1$ with non-negative integer coefficients. 

Set
$\calf_G = \{(u, v) \in (\Zpos)^2 :
u_i(u, v) \in \Zpos,\ i = 0, \ldots, 7\}$.

Note that $p_G \in \calb_{3,1}^+$ only if $a_3 \neq 0$, in which case Definition~\ref{def:frieze} identifies
friezes adapted to $p_G$ with homomorphisms
$f: {\rm ev}_{p_G}(\calu_{3,1}) \to \Z$ satisfying
$f(u_i) > 0$ for all $i$. Such homomorphisms are determined by
their values $(u, v)$ at $(x, y)$, giving
\begin{equation}
    {\rm Frieze}(\calu_{3,1}, p_G)
    \stackrel{\simeq}{\longrightarrow} \calf_G.
\end{equation}

\begin{lemma}\label{lemma:bij-points-gfp}
    The projection map $ Y_G(\Zpos) \longrightarrow (\Zpos)^2$ given by $
    (x, y, z) \mapsto (x, y)$ is a bijection onto $\calf_G$.
\end{lemma}
\begin{proof}
    Five of the eight specialized cluster variables are
    positive integers, being coordinates of the point or quotients
    supplied by the surface equation. The remaining three are handled
    in three steps: Laurent positivity gives nonnegativity, explicit
    polynomial identities give integrality, and a positive power of
    $c$ excludes the value zero.

    Let $(u, v, w) \in Y_G(\Zpos)$. Thus $u_0(u, v) = u$, $u_1(u, v) = v$, and the
    surface equation gives
    \[
        u_2(u, v) = \frac{c + bv}{u}, \qquad
        u_6(u, v) = w, \qquad
        u_7(u, v) =
        \frac{c + a_1u + a_2u^2 + a_3u^3}{v}.
    \]
    These five values are positive integers. The values $u_3(u, v)$, $u_4(u, v)$, and $u_5(u, v)$ are integers because
    \[
        x_4 = x_2x_6 - a_3bx_0 - a_2b,
        \quad
        x_3 = x_2x_4 - a_3bc,
        \quad
        x_5 = x_4x_6 - a_3b^2.
    \]
    On the other hand, the explicit Laurent expressions in
    Example \ref{ex:G2} for $x_3,x_4,x_5$ show $u_3(u,v),u_4(u,v),u_5(u,v)$ are non-negative. To see they are in fact positive, note that these Laurent expressions each contain a positive power of $c$ in the numerator, thus their specialisations $u_3,u_4,u_5$ are 
    non-zero. Periodicity
    now gives
    $u_k(u, v) \in \Zpos$ for every $k$, and hence
    $(u, v) \in \calf_G$.

    Projection to $(u, v)$ is injective because the surface equation
    determines $w$. Conversely, if $(u, v) \in \calf_G$, then
    $u_6(u, v) \in \Zpos$ and the defining formula for $x_6$ gives
    $(u, v, u_6(u, v)) \in Y_G(\Zpos)$.
\end{proof}

The proof of Proposition \ref{prop:finite-case} then reduces to the following lemma, whose method of proof appears already in 
\cite{muller,dSG-Yfrieze}. 
\begin{lemma}\label{lem:gfp-finite}
    For each $G$, the set $\calf_G$ is finite.
\end{lemma}
\begin{proof}
    We multiply the four exchange relations of each parity and
    compare the products of the even and odd entries. Combining
    the two bounds gives an inequality in which the exponent of
    the even product on the left exceeds its exponent on the
    right. Hence both products are bounded.

    Fix $(u, v) \in \calf_G$ and write $u_k = u_k(u, v)$. Set
    $M_0 = u_0u_2u_4u_6$ and $M_1 = u_1u_3u_5u_7$. The exchange
    relations described in Figure \ref{fig:G2-exchange} give
    $M_0^2 = E_1(u_1)E_3(u_3)E_5(u_5)E_7(u_7)$
    and
    $M_1^2 = E_0(u_0)E_2(u_2)E_4(u_4)E_6(u_6)$.
    Set $d = \max(a_3, a_2, a_1, b, c)$. The formulas in
    Figure~\ref{fig:G2-exchange}, together with $u_k \geq 1$, give
    $M_0^2 \leq 2^4d^8M_1 \quad\text{and}\quad M_1^2 \leq
    2^{12}d^{14}M_0^3$.
    Indeed, each exchange polynomial with odd index has two terms, and its
    coefficients are monomials in $a_1, a_2, a_3, b, c$ of degrees at most $1, 3, 3, 1$ for
    $E_1, E_3, E_5, E_7$ respectively; since $u_k \geq 1$, the four
    factors in the first product are bounded by $2du_1$, $2d^3u_3$,
    $2d^3u_5$, and $2du_7$, and $1 + 3 + 3 + 1 = 8$. The second
    estimate is obtained in the same way. Each exchange
    polynomial with even index has four terms and degree $3$ in $u$, and its
    coefficient monomials have degrees at most $1, 3, 6, 4$ for
    $E_0, E_2, E_4, E_6$ respectively. Hence
    $M_1^2 \leq 4^4d^{14}M_0^3 \leq 2^{12}d^{14}M_0^3$, where
    $1 + 3 + 6 + 4 = 14$.

    Squaring the bound for $M_0^2$ and substituting the bound
    for $M_1^2$ gives $M_0^4 \leq 2^{20}d^{30}M_0^3$.
    Since $M_0 > 0$, we obtain $M_0 \leq 2^{20}d^{30}$.
    The bound for $M_1^2$ then gives $M_1 \leq 2^{36}d^{52}$.
    Since $u = u_0 \leq M_0$ and $v = u_1 \leq M_1$,
    only finitely many pairs $(u, v)$ occur.
\end{proof}

\begin{remark}
    \label{rem:effective-finite-type}
    The proof of Lemma~\ref{lem:gfp-finite} is effective: for each
    fixed
    $G$, it bounds every point of $Y_G(\Zpos)$. Mohanty~\cite{mohanty-1977},
    Ramasamy--Mohanty~\cite{ramasamy-mohanty84}, and
    Subburam--Togb\'{e}~\cite{subburam-togbe:G2} obtain sharper
    bounds for several individual equations of bi-degree $(3, 1)$.
\end{remark}

The finiteness in Proposition~\ref{prop:finite-case} is not uniform
in the coefficients.
Indeed, even in bi-degree $(1, 1)$ the number of positive integral points is
unbounded.

\begin{proposition}\label{prop:no-uniform-bound}
    We have
    \[
        \sup_{a,b,c \in \Zpos}
        \#Y_{ax+by+c}(\Zpos)   = \infty.
    \]
\end{proposition}
\begin{proof}
    Set $G_{k} = x + 2^k y + 2^k$ for $k \in \Znn$ and write
    $Y_{k} = Y_{G_{k}}$.
    The coordinate rescaling isomorphism
    \[ 
    Y_k \longrightarrow Y_{k+1}, \quad (x, y, z) \mapsto (2x, y, z)
    \]
    defined in Section \ref{subsec:coord-rescale} restricts to an
    inclusion $Y_k(\Zpos) \hookrightarrow Y_{k + 1}(\Zpos)$. On the
    other hand, the point $(1, 1, 1 + 2^{k + 2}) \in Y_{k + 1}(\Zpos)$
    does not lie in the image of $Y_{k}(\Zpos)$. In particular, the sequence of integers $b_k
    := \#{Y_k(\Zpos)}$ is strictly increasing, which proves our
    claim.
\end{proof}

\begin{remark}
    \label{rem:no-uniform-bound-isomorphism}
    Proposition~\ref{prop:no-uniform-bound} and Example \ref{ex:1-1-single-isom-class} imply that there exists
    no uniform bound on the number of positive integral points of
    Mordell--Schinzel surfaces belonging to the same isomorphism
    class over $\Q$.
\end{remark}

\subsection{Positive constant term: infinite type}\label{ss:infinite-type}
Assume that $m,n \geq 1$. Recall that $\{x_i : i \in \Z\}$ denotes the set of cluster variables of $\calu_{m,n}$. 
\begin{lemdef}[{\cite[\S 3.1]{dhl}}]
\label{lemdef:unitary-frieze}
    Given $p \in \calb_{m,n}^+$ and $i \in \Z$, there exists a unique frieze $f_i: \calu_{m,n} \to \Z$ which is adapted to $p$ and such that 
    \[
        f_i(x_i) = f_i(x_{i+1}) = 1.
    \]
    Friezes arising in this way are called {\it unitary}.
\end{lemdef}

\begin{lemma}\label{lem:rank-two-unitary-separation}
Assume that $mn\geq4$.
    Let $f$ be a unitary frieze adapted to some $p \in \calb_{m,n}^+$, say $f(x_i) = f(x_{i+1}) = 1$ for some fixed $i \in \Z$. For all $j \in
    \Z, j \neq i$ we have $f(x_j) f(x_{j + 1}) > 1$. 
\end{lemma}
\begin{proof}
    To ease notation, set $u_k = f(x_k)$ for all $k \in \Z$. 
    Relabeling if necessary, we may assume that $i = 0$, i.e., $u_0 =
    u_1 = 1$. Suppose there exists $j \in \Z, j \neq 0$ such that
    \begin{equation}\label{eq:unitary-twice}
        u_j\, u_{j + 1} = 1.
    \end{equation}
  
    Applying $f$ to the exchange relation \eqref{eq:cv-recursive} for
    each $k \in \Z$ gives $u_{k - 1}u_{k + 1} = E_k(u_k)$, where $E_k
    \in \Znn[t]$ has degree $d_k \in \{m, n\}$ and has positive
    constant and leading coefficients. Thus,
    \begin{equation}\label{eq:rank-two-growth-ineq}
        u_{k - 1}u_{k + 1} \geq 1 + u_k^{d_k}, \qquad \forall k \in
        \Z.
    \end{equation}
    For convenience, suppose $j > 0$ (the case $j < 0$ is symmetric). Moreover, suppose that $j$ is minimal with
    \eqref{eq:unitary-twice}. Since $u_2 = E_1(1) > 1$, we have $j
    \geq 2$. Let $M = \max\{u_0, u_1, \ldots, u_j, u_{j + 1}\}$.
    This maximum satisfies $M > 1$, and $M = u_\ell$ for some
    $\ell \in [2, j - 1]$. If $d_\ell \geq 2$, then
    \eqref{eq:rank-two-growth-ineq} and the maximality of $M$ give
    \[
        M^2 \geq u_{\ell - 1}u_{\ell + 1}
        \geq 1 + M^{d_\ell} \geq 1 + M^2,
    \]
    a contradiction. Thus $d_\ell = 1$. Set $d_{\ell - 1} = d_{\ell +
    1} = : D$. Since $mn \geq 4$, we have $D \geq 4$.

    If either $u_{\ell - 1} = 1$ or $u_{\ell + 1} = 1$, then $u_{\ell
    - 1}u_{\ell + 1} \leq M$, contradicting $u_{\ell - 1}u_{\ell + 1}
    \geq M + 1$. Thus $u_{\ell - 1}, u_{\ell + 1} > 1$ and
    \eqref{eq:rank-two-growth-ineq}, applied at $\ell - 1$ and $\ell +
    1$, yields $u_{\ell - 1}^D < M^2, \qquad u_{\ell + 1}^D < M^2$.
    Therefore $u_{\ell - 1}u_{\ell + 1} < M^{4/D} \leq M$, again
    contradicting $u_{\ell - 1}u_{\ell + 1} \geq M + 1$.
\end{proof}

\begin{corollary}\label{cor:infinite-frieze}
    Assume that $mn \geq 4$. For each $p \in \calb_{m,n}^+$, the set ${\rm Frieze}(\calu_{m,n}, p)$ is infinite.
\end{corollary}
\begin{proof}
    Lemma \ref{lem:rank-two-unitary-separation} defines an inclusion of sets $\Z \hookrightarrow {\rm Frieze}(\calu_{m,n}, p)$ given by mapping $i$ to the unitary frieze $f_i$ adapted to $p$ constructed in Lemma-Definition \ref{lemdef:unitary-frieze}. 
\end{proof}

\begin{proposition}\label{prop:infinite-case}
    If $mn \geq 4$, then
    $Y_G(\Zpos)$ is infinite.
\end{proposition}
\begin{proof}
    Write $G = A(x) + B(y) -c$ and set $p_G = (a_1,\ldots, a_m,b_1,\ldots, b_n,c) \in \calb_{m,n}^+$ as in Section \ref{sec:friezes}. Theorem \ref{thm:frieze-and-integer-points} defines a bijection from $Y_G(\Zpos)$ to ${\rm Frieze}(\calu_{m,n}, p_G)$, and Corollary \ref{cor:infinite-frieze} shows the latter is infinite. 
\end{proof}

\subsection{Zero constant term}
In this subsection, fix a nonzero pure $G \in \Znn[x, y]$
with $c = 0$, and write $G = A(x) + B(y)$ with
$A(0) = B(0) = 0$ and $A(x) = xA^*(x)$,
$B(y) = yB^*(y)$.

\begin{lemma}
    \label{lemma:degenerate-reducible}
    Suppose that $Y_G$ is not geometrically integral. Then $Y_G(\Zpos)$ is infinite.
\end{lemma}
\begin{proof}
By \eqref{eq:integrality-criterion}, either $x \mid G$ or $y \mid G$, equivalently $B^*=0$ or $A^* = 0$. Suppose $B^* =0$ (the case $A^* =0$ is symmetric), so $G = xA^*$. Since $G \in \Znn[x,y]$, $Y_G(\Zpos)$ contains the infinite set $\{(t,1,A^*(t)) : t \in \Zpos\}$ and we are done.
\end{proof}

We now assume that $Y_G$ is geometrically integral. By \eqref{eq:integrality-criterion}, $AB \neq 0$. As in \S \ref{ss:main-result}, set $ k = \ord_0 A$ and $\ell = \ord_0 B$, so that 
\begin{equation}
    \label{eq:degenerate-orders}
    A(x) = \sum_{i = k}^m a_i x^i,
    \qquad
    B(y) = \sum_{j = \ell}^n b_j y^j,
\end{equation}
where
$a_k a_m b_\ell b_n \neq 0$. Set $A(x) = x^k A_0(x)$ and $B(y) = y^\ell B_0(y)$.

\begin{proposition}
    \label{prop:degenerate-case}
    The set
    $Y_G(\Zpos)$ is finite if and only if
    \[
    \min(m, n) = 1 \quad\text{and}\quad \max(m, n) - \max(k, \ell)
    \leq
    3.
    \]
\end{proposition}

We prove Proposition \ref{prop:degenerate-case} by relating $Y_G$ with $c =0$ to auxiliary Mordell--Schinzel surfaces defined by  polynomials with nonzero constant term, and then applying Propositions \ref{prop:finite-case} and \ref{prop:infinite-case}. We first record these auxiliary surfaces, in the following two lemmas.

\begin{lemma}\label{lem:degen-injections}
    Suppose that $n = 1$, so $B(y) = by$ with $b \neq 0$. There are injections 
\[
\begin{tikzcd}
            Y_{A_0(x)+by}(\Zpos) \arrow[r, hook, "\iota^-"] &  Y_G(\Zpos)  \arrow[r, hook, "\iota^+"] & Y_{b^kA_0(x)+by}(\Zpos)
   \end{tikzcd}
\]
    given by $\iota^-(x, y, z) = (x, x^ky, z)$ and $\iota^+(x,y,z) = (x, b^ky/x^k, z)$.
\end{lemma}
\begin{proof}
    The map $\iota^-$ follows by substituting $y \mapsto x^ky$.
    For the other direction,
    the equation of $Y_G$ forces successive
    divisibility by $x$. Indeed, for $(x,y,z) \in Y_G(\Zpos)$,
    set $y_0=y$ and define $y_{i+1}=by_i/x$. Inductively,
    \[
        y_i\in\Zpos, \qquad
        xy_i z=b^i x^{k-i}A_0(x)+by_i
        \qquad (0\leq i\leq k):
    \]
    for $i<k$, the identity gives $x \mid by_i$, and multiplying 
    it by $b/x$ gives the next identity.
    Thus $y_k=b^ky/x^k$ defines $\iota^+$.
    Both maps are injective since
    their first two coordinates recover $x$ and $y$,
    and they preserve $z$.
\end{proof}
\begin{remark}
    After base change to $\Z[1/b]$, the map $\iota^-$ is an isomorphism.
\end{remark}

\begin{lemma}\label{lem:degen-injection2}
    Suppose that $\min(m,n) \geq 2$ and $k = 1$. 

    \begin{enumerate}
        \item If $\ell \geq 2$, there is an injection 
        \[
            Y_F(\Zpos) \hookrightarrow Y_G(\Zpos)
        \]
        given by $(x,y,z) \mapsto (x^\ell, xy,z + y^\ell C(xy))$, where $C$ is uniquely determined by $B_0(y) = b_\ell + yC(y)$, and 
        \[F = A_0(x^\ell) + b_\ell y^\ell \in \Znn[x,y]
        \]
        has $F(0,0)\neq 0$ and bidegree $(\ell(m-1), \ell)$. 
    \item If $\ell = 1$, there is an injection 
    \[
        Y_H(\Zpos) \hookrightarrow Y_G(\Zpos)
    \]
    given by $(x,y,z) \mapsto (xR(x,y), yR(x,y), z)$ where $R(x,y) = a_1x + b_1y$ and 
    \[
     H(x, y) = 1 +
        \sum_{i = 2}^m a_i x^iR(x, y)^{i - 2} +
        \sum_{j = 2}^n b_j y^jR(x, y)^{j - 2} \in \Znn[x,y]
    \]
    has $H(0,0) = 1$ and bidegree $(2m-2,2n-2)$. 
    \end{enumerate}
\end{lemma}
\begin{proof}
   (1) The substitution $x=u^\ell$, $y=uv$ makes the lowest terms
    $a_1x$ and $b_\ell y^\ell$ divisible
    by the same power $u^\ell$.
    Since $B_0(t)=b_\ell+tC(t)$,
    \[
        A(u^\ell)+B(uv)
        =
        u^\ell\bigl(F(u,v)+uv\,v^\ell C(uv)\bigr).
    \]
    The second term $uv\,v^\ell C(uv)$ is divisible by $uv$
    and is absorbed by translating the third coordinate,
    as in Lemma~\ref{lem:mixed-reduction}.
    Thus, if $(u,v,w)\in Y_F(\Zpos)$, then
    \[
        A(u^\ell)+B(uv)
        =
        (u^\ell)(uv)\bigl(w+v^\ell C(uv)\bigr),
    \]
    which proves that the image of $Y_F(\Zpos)$
    is in $Y_G(\Zpos)$.
    The first coordinate determines $u>0$,
    the second coordinate then determines $v$,
    and the third coordinate determines $w$,
    so the map is injective.
    Finally, $F(0,0)=a_1>0$,
    while $\deg A_0=m-1$,
    which gives the bidegree $(\ell(m-1), \ell)$.

    (2) Write $R=R(u,v)$.
    The choice of $R$ turns the linear part of $G$
    into $R^2$ after the substitution
    $(x,y)=(uR,vR)$:
    \[
        a_1uR+b_1vR=R^2.
    \]
    Consequently,
    \[
        A(uR)+B(vR)=R^2H(u,v).
    \]
    If $(u,v,w)\in Y_H(\Zpos)$, then
    \[
        (uR)(vR)w=R^2uvw=R^2H(u,v),
    \]
    so the map takes $Y_H(\Zpos)$ into $Y_G(\Zpos)$.

    Its image $(x,y,z)$ satisfies
    $a_1x+b_1y=R(u,v)^2$.
    Since $R(u,v)>0$,
    the positive square root recovers $R(u,v)$,
    after which $u=x/R$ and $v=y/R$.
    Hence the map is injective.

    It remains only to verify the bidegree of $H$.
    The polynomial $H$ has nonnegative integral coefficients
    and constant term $1$.
    On $v=0$, its highest-degree term is
    $a_ma_1^{m-2}u^{2m-2}$,
    while on $u=0$ its highest-degree term is
    $b_nb_1^{n-2}v^{2n-2}$.
    Thus $H$ has bidegree $(2m-2,2n-2)$.
\end{proof}

{\it Proof of Proposition \ref{prop:degenerate-case}}.
    We split the proof into four cases, depending on the values of $(m,n,k,\ell)$ associated to $G$.
    
    {\bf Case 1:} $\min(m, n) = 1$. Suppose that $n = 1$ (the case $m = 1$ is symmetric). Then $\ell = 1$, $B(y) = by$ and $A(x) = x^kA_0(x)$ with $b,A_0(0) \neq 0$. Lemma \ref{lem:degen-injections} gives injections $Y_{G_-}(\Zpos) \hookrightarrow Y_G(\Zpos) \hookrightarrow Y_{G_+}(\Zpos)$, 
    where $G_- = A_0(x) + by$ and
    $G_+ = b^kA_0(x) + by$ have positive constant term and
    bi-degree $(m-k,1)$. Propositions \ref{prop:finite-case} and \ref{prop:infinite-case} then imply $Y_G(\Zpos)$ is finite if and only if $m-k = \max(m,n) - \max(k,\ell) \leq 3$. 

    In the remaining three cases, we prove that $Y_G(\Zpos)$ is infinite.
    
    {\bf Case 2:} $\min(m,n) \geq 2$ and $\min(k,\ell) \geq 2$. The points 
    \[
        \left(t,t, \frac{A(t) + B(t)}{t^2}\right), \quad t\in \Zpos
    \]
    are distinct positive integral points of $Y_G$, so $Y_G(\Zpos)$ is infinite.

    {\bf Case 3:} $\min(m,n) \geq 2, \min(k,\ell) = 1$ and $\max(k,\ell) \geq 2$. Suppose that $k = 1$ (the case $\ell = 1$ is symmetric), so $\ell \geq 2$. By Lemma \ref{lem:degen-injection2} there is an injection $Y_F(\Zpos) \hookrightarrow Y_G(\Zpos)$. Moreover, $F(0,0) \neq 0$, $F$ has bidegree $(\ell(m-1), \ell)$ with $\ell^2(m-1) \geq 4$, so Proposition \ref{prop:infinite-case} implies $Y_F$, and thus $Y_G$, have infinitely many integral points.

    {\bf Case 4:} $\min(m,n) \geq 2$ and $ \min(k,\ell) = \max(k,\ell)=  1$. Then $k = \ell = 1$, and Lemma \ref{lem:degen-injection2} implies there exists an injection $Y_H(\Zpos) \hookrightarrow Y_G(\Zpos)$, where $H(0,0) =1$ and $H$ has bidegree $(2m-2,2n-2)$ with $(2m-2)(2n-2) \geq 4$. Applying Lemma~\ref{lem:mixed-reduction} to remove the mixed
    terms of $H$, then Proposition~\ref{prop:infinite-case},
    shows that $Y_H(\Zpos)$, and hence $Y_G(\Zpos)$, is infinite. 
\qed

\begin{proof}[Proof of Theorem~\ref{thm:finiteness-classification}]
    For pure $G$,
    Proposition~\ref{prop:degenerate-case}
    proves part (i),
    and Lemma~\ref{lemma:degenerate-reducible}
    proves part (ii).
    If the pure part of $G$ is nonzero,
    Lemma~\ref{lem:mixed-reduction} gives the general case.
    Otherwise $G = xyH$ with $H \in \Znn[x, y]$ nonzero,
    and $(t, 1, H(t, 1))$, for $t \in \Zpos$, gives
    infinitely many positive integral points.
\end{proof}

\subsection{Zariski density in the indefinite case}
We return to the case when $c > 0$ and
strengthen Proposition~\ref{prop:infinite-case}
to Zariski density of $Y_G(\Zpos)$ when $mn > 4$.
We will use Proposition~\ref{prop:infinite-case}
to produce an infinite sequence of points
$P_k = (x_k, y_k, z_k) \in Y_G(\Zpos)$
such that $x_k \rightarrow \infty$ and
such that the ratios
$\log y_k / \log x_k$ converge to an irrational number.

The following elementary observation then shows that
no algebraic curve can contain infinitely many of these points.
The proof uses a local parameter on the smooth projective model
of the curve, as in the argument of
Corvaja--Zannier~\cite[Section~2]{corvaja-zannier-curves}.

\begin{lemma}
    \label{lem:logarithmic-ratio}
    Let $C$ be an irreducible complex algebraic curve,
    let $(P_k)_{k \in \Zpos}$ be a sequence
    of points in $C(\C)$,
    and let $f$ and $g$ be nonzero rational functions on $C$
    that are defined and nonzero at each point $P_k$.
    If $\verts{f(P_k)} \rightarrow \infty$ and
    $\log \verts{g(P_k)} / \log \verts{f(P_k)}$
    converges to a finite limit $\alpha$ as $k \rightarrow \infty$,
    then $\alpha \in \Q$.
\end{lemma}

\begin{proof}
    Let $\overline{C}$ be a projective completion of $C$
    and let $\nu: \widetilde{C} \to \overline{C}$ be its
    normalization,
    so $\widetilde{C}$ is a smooth projective curve.
    For each $k$, choose $\widetilde{P}_k \in
    \widetilde{C}(\C)$ with $\nu(\widetilde{P}_k) = P_k$.
    View $f$ and $g$ as rational functions
    on $\widetilde{C}$ via $\nu$.
    Since $\widetilde{C}(\C)$ is compact in the
    complex topology,
    a subsequence $\widetilde{P}_{k_i}$
    converges to a point $P \in \widetilde{C}(\C)$.
    The function $f$ has a pole at $P$, since otherwise
    it would be bounded near $P$, contradicting
    $\lvert f(\widetilde{P}_{k_i}) \rvert \rightarrow \infty$.

    Choose a local parameter $t$ at $P$, namely a rational
    function regular at $P$ with a simple zero at $P$.
    We can write $f = t^{-a}u$ and $g = t^{-b}v$
    for some $a \in \Zpos$ and $b \in \Z$, where $u$ and $v$
    are holomorphic and nonzero at $P$.
    Put $t_i = t(\widetilde{P}_{k_i})$, so $t_i \rightarrow 0$
    as $i \rightarrow \infty$.
    Since $u$ and $v$ have nonzero limits along this subsequence,
    we obtain
    \[
        \frac{\log\verts{g(P_{k_i})}}{\log\verts{f(P_{k_i})}}
        = \frac{b\log\verts{t_i}^{-1} + O(1)}
            {a\log\verts{t_i}^{-1} + O(1)}
            \xrightarrow{i \rightarrow \infty}
            \frac{b}{a} \in \Q.
    \]
    This is also the limit $\alpha$ of the original sequence
    $\log\verts{g(P_k)} / \log\verts{f(P_{k})}$,
    so $\alpha \in \Q$.
\end{proof}

We then use Lemma~\ref{lem:logarithmic-ratio}
to show that no algebraic curve in $Y_G$
can contain infinitely many of the points $P_k$ from
Proposition~\ref{prop:infinite-case}.

\begin{proposition}
    \label{prop:indefinite-density}
    Let $G \in \Znn[x, y]$ have bidegree $(m, n)$
    and constant term $c$.
    If $c > 0$ and $mn > 4$, then
    $Y_G(\Zpos)$ is Zariski dense.
\end{proposition}

\begin{proof}
    By Lemma~\ref{lem:mixed-reduction},
    we may assume that $G$ is pure.
    Let $(P_k)_{k \in \Zpos}$ be the infinite sequence
    of pairwise distinct points $P_k = (x_k, y_k, z_k)$
    in $Y_G(\Zpos)$ corresponding to the unitary friezes $f_{2k}$,
    as in the proof of Proposition~\ref{prop:infinite-case}.
    Here $x$ and $y$ are the coordinate functions on $Y_G$,
    so $x_k = x(P_k)$ and $y_k = y(P_k)$.

    First, we prove that $x_k \longrightarrow \infty$
    and that $\log y_k / \log x_k$ converges to the
    irrational number
    \[
        \alpha = \frac{mn - \sqrt{mn(mn - 4)}}{2n};
    \]
    note that $\alpha \in \overline{\Q} - \Q$
    since $mn(mn - 4)$ lies strictly
    between the consecutive squares $(mn - 3)^2$
    and $(mn - 2)^2$ when $mn > 4$.
    Set $\lambda := (m - \alpha)/\alpha > 1$.
    For fixed $k \in \Zpos$,
    write $L_j = \log x_j(P_k)$, with $x_0 = x$ and $x_1 = y$.
    The $L_j$ are nonnegative real numbers such that
    $L_{2k} = L_{2k + 1} = 0$.
    The exchange relations give
    \[
        L_{j - 1} + L_{j + 1} = d_jL_j + O_j(1),
    \]
    where $d_j = m$ for even $j$ and $d_j = n$ for odd $j$.
    The error $O_j(1)$ is positive and bounded independently of $k$.
    Consequently,
    \[
        L_{2j + 3} - \alpha L_{2j + 2}
        > \lambda(L_{2j + 1} - \alpha L_{2j}).
    \]
    Iterating this inequality from $j = 0$ to $j = k - 1$
    and using $L_{2k} = L_{2k + 1} = 0$ gives
    $0 \leq \log y_k < \alpha \log x_k$.
    Since the pairs $(x_k, y_k)$ of positive integers are
    pairwise distinct and satisfy $y_k < x_k^\alpha$,
    we have $x_k \longrightarrow \infty$.
    For each fixed $j \geq 1$,
    solving the recurrence for $L_{2j}$
    and using $L_{2j} \geq 0$ gives
    \[
        0 \leq \alpha - \frac{\log y_k}{\log x_k}
        \leq \frac{m - 2\alpha}{\lambda^{2j} - 1}
            + \frac{C_j}{\log x_k},
    \]
    where $C_j$ is independent of $k$.
    Letting $k \to \infty$ with $j$ fixed and then
    $j \to \infty$ proves that
    $\log y_k / \log x_k \rightarrow \alpha$.
    Thus our infinite sequence of pairwise distinct points
    satisfies $x_k \to \infty$ and
    \[
        \frac{\log y_k}{\log x_k}
        \xrightarrow{k \to \infty}
        \alpha \in \overline{\Q} \setminus \Q.
    \]

    Next, we show that every irreducible curve $C \subset Y_G$
    contains only finitely many of the points $P_k$.
    Suppose that there is an irreducible curve $C \subset Y_G$
    such that $P_k \in C$ for infinitely many $k \in \Zpos$.
    Choose a subsequence $P_{k_i} \in C(\C)$ and set
    $f = \restr{x}{C}$ and $g = \restr{y}{C}$;
    these are regular functions on $C$ that are
    defined and nonzero at every $P_{k_i}$.
    Moreover, $f(P_{k_i}) = x_{k_i} \rightarrow \infty$ and
    \[
        \frac{\log\verts{g(P_{k_i})}}{\log\verts{f(P_{k_i})}}
            = \frac{\log y_{k_i}}{\log x_{k_i}}
                \xrightarrow{i \to \infty} \alpha.
    \]
    Lemma~\ref{lem:logarithmic-ratio} therefore implies
    $\alpha \in \Q$, but this contradicts
    $\alpha \in \overline{\Q} \setminus \Q$.

    Finally,
    a proper closed subset of the integral surface $Y_G$
    has only finitely many curve and point components.
    Therefore the set of points $\{P_k\}_{k \in \Zpos}$
    is Zariski dense in $Y_G$,
    and so is $Y_G(\Zpos)$.
\end{proof}

\section{Orbits of positive integral points}
\label{sec:orbits}
In this section, we prove the orbit and asymptotic counting
assertions of Theorems~\ref{thm:trichotomy}(ii)
and \ref{thm:affine-counting}.
We also prove Zariski density for $G = ax^2 + by^2 + c$
when $a, b, c \in \Zpos$ and $ab > 1$.

For the first two subsections, fix a pure $G \in \Znn[x,y]$
with decomposition \eqref{eq:pure-poly} and satisfying the conditions $mn = 4$, $c \neq 0$ and $abc^{m+n-4} = 1$. The last condition is equivalent to assuming either \eqref{eq:am-bn-c=1} or \eqref{eq:m-n=2-a-b=1} holds. There are two possibilities for $G$ (up to swapping $x$ and $y$):
\begin{align}\label{eq:bi-quadratic}
    G &= x^2 + a_1x + y^2 + b_1y + c, \qquad &a_1,b_1 \in \Znn, c \in \Zpos \\
    G&= x^4 + a_3x^3 + a_2x^2 + a_1x + y + 1, \quad &a_1,a_2,a_3 \in \Znn. \label{eq:quartic-linear}
\end{align}
We refer to \eqref{eq:bi-quadratic} as the {\it bi-quadratic case}, and to \eqref{eq:quartic-linear} as the {\it quartic-linear case}. Our proofs of Theorem~\ref{thm:trichotomy}(ii) and Theorem~\ref{thm:affine-counting} treat these two cases separately.

\subsection{The bi-quadratic case}
\label{sec:quadratic-case}
Suppose that $G$ is of the form \eqref{eq:bi-quadratic}. In the notation of \S \ref{subsec:companion}, we have that $A^\dagger = A$ and $B^\dagger = B$, so $\calg_{A,B}$ consists of the single object $Y_{G}$ and $\Gamma_{A,B}$ is generated by $\mu_{A,B,1}$ and $\mu_{A,B,2}$. For simplicity, we write $\mu_{A,B,i} = \mu_i, i = 1,2$.

\begin{lemma}
    \label{lem:quadratic-mutations}
    For $i = 1,2$, we have  
    \[
        \mu_i^2 = \Id, \quad \mu_i(Y_G(\Zpos)) = Y_G(\Zpos), \quad \mu_i^*(z) = z,
    \]
    where $\mu_i^*$ denotes the pullback of the automorphism $\mu_i$. 
\end{lemma}
\begin{proof}
   Substituting the formulas for $A$ and $B$ in \eqref{eq:mu_1} and \eqref{eq:mu_2}, we get
   \[
    \mu_1(x,y,z) = (yz - x - a_1, y,z) \quad \text{and} \quad \mu_2(x,y,z) = (x,xz - y - b_1,z).
   \]
   Note that $yz - x- a_1 = \frac{y^2 + b_1y + c}{x}$ and $xz - y - b_1 = \frac{x^2 + a_1x + c}{y}$. The claims follow.
\end{proof}

Set $\tau = \mu_2 \circ \mu_1 \in \Gamma_{A,B}$. Each $P = (x_0,y_0,z_0) \in Y_G(\Zpos)$ defines a bi-infinite sequence ${\bf u}(P) := (u_k)_{k \in \Z}$ via 
\[
    (u_{2k}, u_{2k+1}, z_0) = \tau^k(x_0,y_0,z_0), \quad k \in \Z.
\]
Using the formulas contained in the proof of Lemma \ref{lem:quadratic-mutations}, we have for all $k \in \Z$ 
\begin{equation}\label{eq:mutation-2-2}
        u_ku_{k + 2} =
    \begin{cases}
        u_{k + 1}^2 + b_1u_{k + 1} + c,
        &k\text{ even}, \\
        u_{k + 1}^2 + a_1u_{k + 1} + c,
        &k\text{ odd}.
    \end{cases}
\end{equation}

\begin{proposition}
    \label{prop:quadratic-affine-orbits}
    Suppose that $G$ is of the form \eqref{eq:bi-quadratic}. The following statements hold.
    \begin{enumerate}
        \item $Y_{A,B}(\Zpos)$ decomposes into only finitely many $\Gamma_{A,B}$-orbits.
        \item Every $\Gamma_{A,B}$-orbit of $Y_{A,B}(\Zpos)$ is infinite.
        \item Every $\Gamma_{A,B}$-orbit of $Y_{A,B}(\Zpos)$ is contained in a level curve of $z$.
    \end{enumerate}
\end{proposition}

\begin{proof}
The final assertion follows from
    Lemma \ref{lem:quadratic-mutations}. Since $z$ is
    nonconstant, each level set is at most one-dimensional.

    Our proof of (1) is similar to that in \cite{mills}. Fix $P \in Y_{A,B}(\Zpos)$ and ${\bf u}(P) = (u_k)_{k \in \Z}$. Let $j \in \Z$ such that $u_j = \min(u_k : k \in \Z)$, and write $u = u_j$. Then
    $u_{j - 1} = u + \delta_- \quad\text{and}\quad u_{j + 1} = u +
    \delta_+$
    for some $\delta_-, \delta_+ \in \Znn$. We first show that $\delta_{\pm}$ are bounded above  by a constant that is {\it independent of }${\bf u}(P)$. Set $d = b_1$ if $j$ is odd and $d = a_1$ if $j$ is even. Then \eqref{eq:mutation-2-2} implies
    \[
        u_{j-1} u_{j+1}=  (u + \delta_-)(u + \delta_+)=  u^2 + du + c.
    \]
    Set $r =
    \delta_- + \delta_+$ and $s = \delta_-\delta_+$. Re-arranging gives
    \begin{equation}
        \label{eq:quadratic-minimum}
        (r - d)u + s - c = 0.
    \end{equation}
    If $r \leq d$, then $\delta_\pm \leq d$, and we are done. Thus assume $r > d$. But then $(r - d)u = c
    - s \leq c \leq cu$, and so $r \leq d + c$. Therefore $\delta_{\pm} \leq d + c$ in all cases. The choice of $d$ depends on ${\bf u}(P)$, but $d \leq a_1+b_1$ always holds: thus there are at most $(a_1+b_1+c+1)^2$ pairs $(\delta_-, \delta_+)$ that can occur.

    Fix an arbitrary $(\delta_-, \delta_+) \in [0,a_1+b_1+c+1]^2$. We now show that there are only finitely many values of $u$ with such a pair. First note that \eqref{eq:quadratic-minimum} determines
    $u$, unless $r = d$ and $s = c$. In that exceptional case, let $e$
    be such that $d+e = a_1+b_1$. Then \eqref{eq:mutation-2-2} gives $uu_{j +
    2} = (u + \delta_+)^2 + e(u + \delta_+) + c$. Since $u >0$, $u
    \mid \delta_+^2 + e\delta_+ + c$, so $u$ is bounded. All in all, we have shown there are only finitely many triples $(u,\delta_+,\delta_-)$. Since $u_{j+1}= u + \delta_+$, there are only finitely many pairs $(u_j,u_{j+1})$, and so only finitely many $\Gamma_{A,B}$-orbits.

    To prove (2), suppose for a contradiction that there exists $P \in Y_G(\Zpos)$ whose $\Gamma_{A,B}$-orbit is finite. Then ${\bf u}(P) = (u_j)_{j\in \Z}$ is periodic. Set 
    \[
        M_0 = \max(u_{2j} : j \in \Z) \quad M_1 = \max(u_{2j+1} : j \in \Z),
    \]
    and choose $i,j \in \Z$ such that $u_{2i} = M_0$ and $u_{2j+1} = M_1$. Then, 
    \[
        M_1^2 \geq u_{2i-1}u_{2i+1} = (u_{2i})^2 + a_1u_{2i} + c > M_0^2.
    \]
    A similar calculation shows $M_0^2 > M_1^2$, which is a contradiction.
 \end{proof}

Retain the height $h(x,y,z) = \max(x,y)$ defined in \S \ref{ss:main-result} and fix a $\Gamma_{A,B}$-orbit $\calo$ in $Y_{G}(\Zpos)$. Pick $P \in \calo$ and set ${\bf u}(P) = (u_{k})_{k\in \Z}$. Set 
\[
P_k= (w_k,y_k,z_\calo), \quad q_k = (w_{k+1},y_k,z_\calo), \quad w_k = u_{2k}, \quad y_k = u_{2k+1},
\]
where $z_\calo =z(P)$ denotes the common $z$-coordinate on $\calo$. Note that $\{P_k, q_k\}_{k \in \Z}$ exhaust $\calo$. 

Consider the graph $T_\calo$ whose vertices are the points of $\calo$, and where two vertices $t_1$ and $t_2$ are joined by an edge if 
\[
\mu_i(t_1) = t_2 \text{ for some } i\in \{1,2\}, \text{ and } t_1 \neq t_2.
\] 
Thus $T_\calo$ is connected, infinite, has no loops, and has each of its vertices of degree at most two. In particular, $T_\calo$ is either a 2-regular tree or a ray; set $e_\calo = 2$ or $1$, respectively. Note both cases occur, as Example \ref{ex:fixed-points} illustrates.

\begin{example}\label{ex:fixed-points}
    \begin{enumerate}
        \item If $G = x^2 + y^2 + 1$, then $\Gamma_{A,B} = \langle \mu_1,\mu_2\rangle$ where 
    \[
        \mu_1(x,y,z) = (yz-x,y,z) \quad \mu_2(x,y,z) = (x,xz-y,z).
    \]
    There is only one $\Gamma_{A,B}$-orbit $\calo$, and $\mu_i$ has no fixed points on $\calo$. Pick $P_0 = (1,1,3) \in \calo$. Then $T_\calo$ is the 2-regular tree
    \[
        \begin{tikzcd}
            \cdots \arrow[r, no head, "\mu_1"] & q_{-1}\arrow[r, no head, "\mu_2"]& P_0 \arrow[r, no head, "\mu_1"] & q_0\arrow[r, no head, "\mu_2"] & P_1 \arrow[r, no head, "\mu_1"] &q_1\arrow[r, no head, "\mu_2"] & \cdots
        \end{tikzcd}
    \]
        \item Take $G = x^2 + y^2 + 3$. Then $\Gamma_{A,B} = \langle \mu_1,\mu_2\rangle$ where $\mu_1$ and $\mu_2$ are given by the same equations as above. Note that $\mu_2(1,2,4) = (1,2,4)$. Thus, the graph $T_\calo$ of the $\Gamma_{A,B}$-orbit $\calo$ of $P_0 = (1,2,4) \in Y_G(\Zpos)$ is the ray
    \[
        \begin{tikzcd}
            P_0 \arrow[r, no head, "\mu_1"] & q_0\arrow[r, no head, "\mu_2"] & P_1 \arrow[r, no head, "\mu_1"] &q_1\arrow[r, no head, "\mu_2"] & \cdots
        \end{tikzcd}
    \]
    \end{enumerate}
\end{example}

\begin{proposition}\label{prop:orbit-count-2-2}
    The coordinates of $P_k$ and $q_k$ satisfy the linear recurrence relation
    \begin{equation}\label{eq:linear-rec-rel-2-2}
        \begin{aligned}
            w_{j + 2} &= (z_{\calo}^2 - 2)w_{j + 1} - w_j
                - (2a_1 + b_1z_{\calo}),\\
            y_{j+2} &= (z_{\calo}^2 - 2)y_{j + 1} - y_j
                - (2b_1 + a_1z_{\calo}).
        \end{aligned}
    \end{equation}
    In particular,
        $\#\{Q\in\calo:h(Q)\leq N\}
        =\alpha_{\calo}\log N+O_{\calo}(1)$,
    where
        \[
            \alpha_{\calo}=\frac{2e_{\calo}}{\log\lambda_{\calo}},
            \qquad
            \lambda_{\calo}=\frac{z_{\calo}^2-2+
                z_{\calo}\sqrt{z_{\calo}^2-4}}{2};
        \]    
\end{proposition}
\begin{proof}
    Note that
    \[
        w_{j+1}+w_j=z_{\calo}y_j-a_1,\qquad
        y_{j+1}+y_j=z_{\calo}w_{j+1}-b_1.
    \]
    Eliminating $y_j$ and $w_{j+1}$ gives
    \eqref{eq:linear-rec-rel-2-2}.

    The surface equation gives
    \[
        (z_{\calo}-2)xy=(x-y)^2+a_1x+b_1y+c>0,
    \]
    so $z_{\calo}\geq3$. In turn this implies the characteristic polynomial of \eqref{eq:linear-rec-rel-2-2} has roots $\lambda_\calo^{\pm 1}$ with $\lambda_\calo >1$ and $0<\lambda_\calo^{-1}<1$. Solving \eqref{eq:linear-rec-rel-2-2} explicitly, we have
    \[
        w_j=C_+\lambda_{\calo}^j+C_-\lambda_{\calo}^{-j}+O_{\calo}(1),
        \qquad
        y_j=D_+\lambda_{\calo}^j+D_-\lambda_{\calo}^{-j}+O_{\calo}(1),
    \]
    for some $C_\pm, D_\pm$. We claim that $C_\pm, D_\pm$ are positive. Suppose not, say $C_+$ is non-positive. If $C_+ <0$, then $w_j <0$ for large enough $j$. If $C_+ = 0$, then for large enough $j$, the sequence $(w_{j})$ is constant. This in turn would imply that the $y_j$ are eventually constant. But then \eqref{eq:linear-rec-rel-2-2} implies the whole orbit is finite,
    contrary to Proposition~\ref{prop:quadratic-affine-orbits}. A similar argument holds for $C_-, D_\pm$. 
    Therefore
    \[
        \log h(P_j)=|j|\log\lambda_{\calo}+O_{\calo}(1),\qquad
        \log h(q_j)=|j|\log\lambda_{\calo}+O_{\calo}(1).
    \]
    Counting the integers $j$ with $h(P_j) \leq N$ and the
    integers $j$ with $h(q_j) \leq N$ separately gives
    \[
        \#\set{j \in \Z : h(P_j) \leq N}
        + \#\set{j \in \Z : h(q_j) \leq N}
        = \frac{4}{\log\lambda_{\calo}}\log N + O_{\calo}(1).
    \]
    If $T_\calo$ is a 2-regular tree, all occurrences in $(P_j)$
    and $(q_j)$ represent distinct points. If $T_\calo$ is a
    ray, every point other than the endpoint occurs twice.
    Therefore
    $\alpha_{\calo} = 2e_{\calo}/\log\lambda_{\calo}$.
\end{proof}

\subsection{The quartic--linear case}
\label{sec:quartic-case}
Suppose that $G$ is of the form \eqref{eq:quartic-linear}. In the
notation of \S\ref{subsec:companion}, we have
\[
    A(x)=x^4+a_3x^3+a_2x^2+a_1x+1,\qquad B(y)=y+1,
\]
with $A^\dagger(x)=x^4A(1/x)$ and $B^\dagger=B$. Thus
$\calg_{A,B}$ has one object if $a_1=a_3$ and two otherwise.
Write $Y=Y_{A,B}$, $Y^\dagger=Y_{A^\dagger,B}$, and
\[
    \mu_i=\mu_{A,B,i},\qquad
    \mu_i^\dagger=\mu_{A^\dagger,B,i}\quad(i=1,2).
\]
The groupoid is generated by
\[
\begin{tikzcd}
    Y \arrow[l,loop left, "\mu_{2}"] \arrow[r, shift left, "\mu_{1}"]& Y^\dagger \arrow[l,shift left, "\mu_1^\dagger"] \arrow[r,loop right, "\mu_2^\dagger"]
\end{tikzcd}
\]
When $A=A^\dagger$, the two objects and their corresponding maps
are identified.

Given $\sigma$ a morphism of $\calg_{A,B}$, its domain and codomain are called the {\it source} and {\it target} of $\sigma$, respectively. 
\begin{lemma}\label{lem:quartic-companion-mutations}
    We have
    \[
        \mu_2^2=\Id,\qquad (\mu_2^\dagger)^2=\Id,\qquad
        \mu_1^{-1}=\mu_1^\dagger.
    \]
    Each map induces a bijection between the positive integral
    points of its source and target.
\end{lemma}
\begin{proof}
    Apply Proposition~\ref{prop:integral-normalization}, using
    $(A^\dagger)^\dagger=A$ and $B^\dagger=B$.
\end{proof}

The vertex group of $Y$ is
\begin{equation}\label{eq:quartic-vertex-group}
    \Gamma_{A,B} = \begin{cases}
        \langle \mu_1,\mu_2\rangle &\text{if } A = A^\dagger, \\
        \langle \mu_2, \mu_1^\dagger \circ \mu_2^\dagger \circ \mu_1\rangle &\text{if } A^\dagger \neq A,
    \end{cases}
\end{equation}
For $P\in Y(\Zpos)$, the groupoid orbit
$\calg_{A,B}\cdot P$ is the set of images of $P$ under
morphisms with source $Y$. It satisfies
\begin{equation}\label{eq:companion-orbit-intersection}
    (\calg_{A,B}\cdot P)\cap Y=\Gamma_{A,B}\cdot P \quad \text{and} \quad (\calg_{A,B}\cdot P)\cap Y^\dagger=\Gamma_{A^\dagger,B}\cdot \mu_1(P).
\end{equation}
The map $\mu_1$ identifies the first intersection with the second. Hence a vertex-group
orbit is finite if and only if its groupoid orbit is finite.

Set $\tau=\mu_2^\dagger\circ\mu_1$ and
$\tau^\dagger=\mu_2\circ\mu_1^\dagger$. As in
\S\ref{sec:quadratic-case}, each $P\in Y(\Zpos)$ defines a
bi-infinite sequence ${\bf v}(P)=(v_k)_{k\in\Z}$, now by
\[
    P_{2k}=(\tau^\dagger\circ\tau)^k(P),\qquad
    P_{2k+1}=\tau(P_{2k}),\qquad
    P_j=(v_{2j},v_{2j+1},*) .
\]
Here, and below, the third coordinate is determined by the first two
and the surface equation. Set $A_j=A$ for even $j$ and
$A_j=A^\dagger$ for odd $j$. Then $P_j\in Y_{A_j,B}$ and
\[
    q_j:=\mu_{A_j,B,1}(P_j)=(v_{2j+2},v_{2j+1},*)\in Y_{A_{j+1},B}.
\]
The points $P_j,q_j$ exhaust $\calg_{A,B}\cdot P$, and
\begin{equation}\label{eq:quartic-orbit-recurrence}
    v_{2j}v_{2j+2}=v_{2j+1}+1,\qquad
    v_{2j-1}v_{2j+1}=A_j(v_{2j}),\qquad j\in\Z.
\end{equation}

Define
\begin{equation}
    \label{eq:invariant-w-4-1}
    \Omega =
    yz^2
    -(x^3 + a_3 x^2 + a_2 x + a_1)z
    -(x^2 + a_3 x + a_2)
    \in \Z[Y],
\end{equation}
and define $\Omega^\dagger \in \Z[Y^\dagger]$ by interchanging $a_3$ and $a_1$.  It is clear that $\Omega$ and $\Omega^\dagger$ are nonconstant.

\begin{lemma}
    \label{lem:quartic-companion-invariant}
We have
    \[
        \mu_2^*\Omega=\Omega, \quad (\mu_2^\dagger)^*\Omega^\dagger = \Omega^\dagger, \quad 
        \mu_1^*\Omega^\dagger=\Omega, \quad (\mu_1^\dagger)^*\Omega = \Omega^\dagger.
    \]
\end{lemma}
\begin{proof}
On $D(xy)\subset Y$, with $x'=(y+1)/x$ and $y'=A(x)/y$,
    the surface equation gives
    \begin{equation}\label{eq:quartic-invariant-forms}
    \begin{aligned}
        \Omega
        &=\frac{(y+1)^2+a_1x(y+1)+a_2x^2+a_3x^3+x^4}{x^2y}\\
        &=\frac{(x')^2+a_1x'+a_2+a_3x+x^2}{y}
         =\frac{y+y'+2+a_1x}{x^2}.
    \end{aligned}
    \end{equation}
    The last expression is unchanged by $\mu_2$, which interchanges
    $y,y'$. The middle expression is unchanged by simultaneously
    interchanging $x,x'$ and $a_1,a_3$.  These identities extend to $Y$ since $D(xy)$ is dense in $Y$. This proves the first and third identities in the statement of the Lemma. The arguments for the second and fourth are symmetric, and we omit them.
\end{proof}

\begin{proposition}\label{prop:quartic-affine-orbits}
    Suppose that $G$ is of the form \eqref{eq:quartic-linear}.
    \begin{enumerate}
        \item $Y_{A,B}(\Zpos)$ decomposes into only finitely many
        $\Gamma_{A,B}$-orbits.
        \item Every $\Gamma_{A,B}$-orbit in $Y_{A,B}(\Zpos)$ is infinite.
        \item Every such orbit is contained in a level curve of $\Omega$.
    \end{enumerate}
\end{proposition}
\begin{proof}
    The third assertion is Lemma \ref{lem:quartic-companion-invariant}.
    For (1), fix $P \in Y(\Zpos)$ and choose $j$ such that
    $x=v_{2j}$ is minimal among the
    even entries of ${\bf v}(P)$, and set
    \[
        v_{2j-2}=x+\delta_-,\qquad v_{2j+2}=x+\delta_+,
        \qquad r=\delta_-+\delta_+,\quad t=\delta_-\delta_+.
    \]
    Here $\delta_\pm\in\Znn$. Write
    $A_j(T)=T^4+\rho T^3+a_2T^2+\sigma T+1$, where
    $(\rho,\sigma)\in\{(a_3,a_1),(a_1,a_3)\}$.
    \eqref{eq:quartic-orbit-recurrence} gives
    \[
        \bigl(x(x+\delta_-)-1\bigr)
        \bigl(x(x+\delta_+)-1\bigr)=A_j(x),
    \]
    hence
    \begin{equation}\label{eq:quartic-minimum-polynomial}
        (r-\rho)x^2+(t-a_2-2)x-(r+\sigma)=0.
    \end{equation}
    If $x=1$, then $t=\rho+a_2+\sigma+2>0$, so only finitely many
    pairs $(\delta_-,\delta_+)$ occur. If $x\geq2$, then
    \[
        r\leq\frac{\rho x^2+(a_2+2)x+\sigma}{x^2-1}
        \leq\frac{4\rho+2a_2+4+\sigma}{3},
    \]
    again leaving only finitely many pairs. For any fixed pair, the
    polynomial in \eqref{eq:quartic-minimum-polynomial} is nonzero:
    otherwise $r=\rho=-\sigma=0$, whence $t=0$, contradicting
    $t=a_2+2$. Thus $x$ also has only finitely many possibilities.
    The polynomial $A_j$ and the pair
    $(v_{2j},v_{2j+1})=(x,x(x+\delta_+)-1)$ determine the groupoid
    orbit. Equation~\eqref{eq:companion-orbit-intersection} proves (1).

    For (2), the identity \eqref{eq:companion-orbit-intersection} implies that a vertex-group orbit is finite if and only if the corresponding
    groupoid orbit is finite. We show that $\calg_{A,B}\cdot P$ cannot be finite. Suppose, for a contradiction, that it is. The bi-infinite sequence ${\bf v}(P) = (v_k)_{k \in \Z}$ must then be periodic. Set $M_0, M_1$ be the maxima of its even and odd
    sequence entries. Applying \eqref{eq:quartic-orbit-recurrence}
    at entries attaining these maxima gives
    \[
        M_0^2\geq M_1+1>M_1,\qquad
        M_1^2\geq A_j(M_0)>M_0^4,
    \]
    a contradiction.
\end{proof}

Fix a $\Gamma_{A,B}$-orbit $\calo$ in $Y(\Zpos)$ and choose
$P\in\calo$. Retain $P_j,q_j,A_j$ from above and set
$w_j=v_{2j}$, $y_j=v_{2j+1}$. Then
\begin{equation}\label{eq:quartic-companion-relations}
    w_jw_{j+1}=y_j+1,\qquad y_{j-1}y_j=A_j(w_j),
\end{equation}
and
\[
    \calo=\{P_j,q_{j-1}:j\in\Z,\ A_j=A\}.
\]
Let $\omega_{\calo}=\Omega(P)$; by
Lemma~\ref{lem:quartic-companion-invariant}, this is also the value
of $\Omega^\dagger$ on the orbit in $Y^\dagger$.
Define $T_{\calo}$ and $e_{\calo}$ as in
\S\ref{sec:quadratic-case}, but with $\mu_1$ and $\mu_2$ replaced by the generators of $\Gamma_{A,B}$ given
in \eqref{eq:quartic-vertex-group}. For the same reasons as in \S\ref{sec:quadratic-case}, $T_{\calo}$ is either a 2-regular tree or a ray, and $e_{\calo}$ is either 2 or 1, respectively.

\begin{proposition}\label{prop:orbit-count-4-1}
    The coordinates $w_j$ satisfy
    \begin{equation}\label{eq:alternating-recurrence-4-1}
        w_{j+1}=\omega_{\calo}w_j-w_{j-1}-\ell_j,
    \end{equation}
    where $\ell_j$ is the coefficient of $T$ in $A_j$. In particular,
    \begin{equation}\label{eq:linear-rec-rel-4-1}
        w_{j+2}=(\omega_{\calo}^2-2)w_j-w_{j-2}
        -\omega_{\calo}\ell_j-2\ell_{j+1},
    \end{equation}
    and
    \[
        \#\{Q\in\calo:h(Q)\leq N\}
        =\alpha_{\calo}\log N+O_{\calo}(1),
    \]
    where
    \[
        \alpha_{\calo}=\frac{\varepsilon_Ae_{\calo}}
            {2\log\lambda_{\calo}},\qquad
        \lambda_{\calo}=\frac{\omega_{\calo}+
            \sqrt{\omega_{\calo}^2-4}}{2},
        \qquad
        \varepsilon_A=\begin{cases}2,&A=A^\dagger,\\1,&A\ne A^\dagger.
        \end{cases}
    \]
\end{proposition}
\begin{proof}
    Equations~\eqref{eq:quartic-invariant-forms} and
    \eqref{eq:quartic-companion-relations} give
    \[
        \omega_{\calo}w_j^2=y_{j-1}+y_j+2+\ell_jw_j
        =w_j(w_{j-1}+w_{j+1}+\ell_j).
    \]
    Division by $w_j$ proves \eqref{eq:alternating-recurrence-4-1}.
    Applying it twice and using $\ell_{j-1}=\ell_{j+1}$ gives
    \eqref{eq:linear-rec-rel-4-1}.

    We follow the proof of Proposition~\ref{prop:orbit-count-2-2}.
    First, $\omega_{\calo}$ is a positive integer by
    \eqref{eq:invariant-w-4-1} and \eqref{eq:quartic-invariant-forms}.
    If $\omega_{\calo}\leq2$, then
    \[
        w_{j+1}-2w_j+w_{j-1}
        =(\omega_{\calo}-2)w_j-\ell_j\leq0.
    \]
    A positive concave sequence indexed by $\Z$ is constant, since
    a nonzero successive difference forces negative entries in one
    direction. Then $y_j=w_jw_{j+1}-1$ is constant as well,
    contradicting Proposition~\ref{prop:quartic-affine-orbits}.
    Thus $\omega_{\calo}\geq3$.

    The recurrence \eqref{eq:alternating-recurrence-4-1} has the
    period-two particular solution
    $r_j=(\omega_{\calo}\ell_j+2\ell_{j+1})/(\omega_{\calo}^2-4)$.
    Hence
    \[
        w_j=r_j+C_+\lambda_{\calo}^j+C_-\lambda_{\calo}^{-j},
        \qquad C_+,C_->0.
    \]
    Positivity follows as in Proposition~\ref{prop:orbit-count-2-2}:
    a negative coefficient contradicts positivity of $w_j$; if one
    vanished, the two integral parity subsequences would converge
    and become constant, forcing a finite groupoid orbit.
    Since $y_j=w_jw_{j+1}-1$,
    \[
        \log h(P_j)=2|j|\log\lambda_{\calo}+O_{\calo}(1),\qquad
        \log h(q_{j-1})=2|j|\log\lambda_{\calo}+O_{\calo}(1).
    \]
    The condition $A_j=A$ allows all integers $j$ if $A=A^\dagger$
    and only even $j$ otherwise. Therefore
    \[
    \begin{aligned}
        &\#\{j:A_j=A,\ h(P_j)\leq N\}
        +\#\{j:A_j=A,\ h(q_{j-1})\leq N\}\\
        &\qquad=\frac{\varepsilon_A}{\log\lambda_{\calo}}\log N+O_{\calo}(1).
    \end{aligned}
    \]
    The repetition count from \S\ref{sec:quadratic-case} applies to
    $T_{\calo}$: multiplying by $e_{\calo}/2$ proves the formula.
\end{proof}

\begin{proof}[Proof of Theorem~\ref{thm:trichotomy}(ii) and Theorem~\ref{thm:affine-counting}]
   The infinitude assertion of Theorem~\ref{thm:trichotomy}(ii) follows from Proposition~\ref{prop:infinite-case}. Under $abc^{m+n-4}=1$, Propositions~\ref{prop:quadratic-affine-orbits} and~\ref{prop:quartic-affine-orbits} give finitely many $\Gamma_{A,B}$-orbits, each contained in a curve. Since $\Gamma_{A,B}\subseteq\Aut_\Z^+(Y_G)$, this proves the remaining assertions of Theorem~\ref{thm:trichotomy}(ii) for pure $G$.

Summing the estimates in Propositions~\ref{prop:orbit-count-2-2} and~\ref{prop:orbit-count-4-1} over these finitely many orbits gives Theorem~\ref{thm:affine-counting}, with $C_G=\sum_{\calo}\alpha_{\calo}>0$. The sum is nonempty since $Y_G(\Zpos)$ is nonempty. Lemma~\ref{lem:mixed-reduction} extends both conclusions to general $G$: the isomorphism preserves positive integral points and their heights, and carries curves to curves.
\end{proof}

\subsection{Nonmonic bi-quadratic surfaces}\label{sec:sharpness}
In this subsection, we fix $G = ax^2 + by^2 + c$ with $a,b,c \in \Zpos$. 

We begin by stating a special case of \cite[Lemma 15]{schinzel-2015}. For completeness, we furnish a proof.
\begin{lemma}
    \label{lem:pell-equation}
    Fix $z \in \Zpos$ and suppose that $D := z^2 - 4ab >0$ is not a square. If there exists a pair $(x,y) \in (\Zpos)^2$ such that $ax^2 - zxy + by^2 = -c$, there are infinitely many such pairs.
\end{lemma}
\begin{proof}
    Set $u = 2ax - zy$. The equation becomes $u^2 - Dy^2 = -4ac$. Let
    $\epsilon > 1$ be a totally positive unit of norm $1$ in the real
    quadratic order $\Z[\sqrt D]$. Since the unit group of the finite
    ring $\Z[\sqrt D]/2a\Z[\sqrt D]$ is finite, after replacing
    $\epsilon$ by a positive power we may assume that $\epsilon \equiv
    1 \pmod{2a\Z[\sqrt D]}$. If $\xi = u + y\sqrt D$ is the element
    associated to the given solution, write $\xi\epsilon^r = u_r +
    y_r\sqrt D$. The elements $u_r$ and $y_r$ are integral, the norm
    equation is preserved, and
    $u_r + zy_r \equiv u + zy \equiv 0 \pmod{2a}$. Thus
    $x_r = (u_r + zy_r)/(2a)$ is an integer. The original quadratic
    inequality places $x/y$ strictly between the two roots of $aT^2 -
    zT + b$. Hence $\xi = 2ax - y(z - \sqrt D) > 0$, while its
    conjugate is negative. Both embeddings of $\epsilon$ are positive.
    It follows that $y_r > 0$, and $u_r + zy_r = (u_r + y_r\sqrt D) +
    (z - \sqrt D)y_r > 0$. Thus $x_r > 0$. The resulting solutions are
    distinct.
\end{proof}

The following lemma amounts to a more precise description of the unitary friezes of $\calu_{2,2}$ adapted to $Y_G$; see \S \ref{ss:infinite-type}. Its proof is similar to that of \cite[Lemma 14]{schinzel-2015}. 

\begin{lemma}
    \label{lem:nonmonic-fiber-points}
    Fix $M = ab$, and suppose that $M > 1$. For every $k \geq 1$, there exists a triple $P_k = (x_k,y_k,z_k)$ such that i) $P_k \in Y_G(\Zpos)$, and ii) $ z_k = 1 + M + aM^{k - 1}c$.
\end{lemma}

\begin{proof}
Set $C_j=aM^{j-1}c$ and
\[
 \widetilde Y_j=(uvw=u^2+Mv^2+C_j),\qquad j\geq1.
\]
In Example~\ref{ex:2-2-isoms}, the identity for odd $j$ and
$(x,y,z)\mapsto(y,x,z)$ for even $j$ identify $Y_{g_j}$ with
$\widetilde Y_j$. The successive isomorphisms become
\[
 (u,v,w)\longmapsto(uw-Mv,u,w):
 \widetilde Y_j\longrightarrow\widetilde Y_{j+1}.
\]
Fix $k\geq1$ and start with
$(u_k,v_k,z_k)=(1,1,1+M+C_k)\in\widetilde Y_k(\Zpos)$.
Apply the inverse maps recursively:
\[
 u_j=v_{j+1},\qquad
 v_j=\frac{v_{j+1}z_k-u_{j+1}}{M},
 \qquad j=k-1,\ldots,1.
\]
The equation of $\widetilde Y_{j+1}$ gives
\[
 u_{j+1}(v_{j+1}z_k-u_{j+1})
   =M(v_{j+1}^2+C_j).
\]
If $u_{j+1}$ is prime to $M$, then $v_j$ is an integer and
\[
 v_j=\frac{v_{j+1}^2+C_j}{u_{j+1}}>0.
\]
For $j\geq2$, we have $M\mid C_j$, so reduction modulo each
prime divisor of $M$ shows that $u_jv_j$ remains prime to $M$.
Starting from $u_kv_k=1$, this proves that the recursion reaches
$\widetilde Y_1(\Zpos)$. At the last step, the same argument
modulo the prime divisors of $a$, using $a\mid C_1$, gives
$\gcd(u_1v_1,a)=1$; this also holds when $k=1$.

Finally,
\[
 u_1(v_1z_k-u_1)=a(bv_1^2+c).
\]
Thus $x=(v_1z_k-u_1)/a$ is a positive integer, and
$(x,v_1,z_k)\in Y_G(\Zpos)$ by the inverse of
$\mu_{g_0,1}$ in Example~\ref{ex:2-2-isoms}.
\end{proof}

\begin{proposition}
    \label{prop:nonmonic-quadratic-density}
    Let $a, b, c \in \Zpos$, with $ab > 1$. There are infinitely
    many fibers of
    $z: Y_G \longrightarrow \A^1$ that contain
    infinitely many positive integral points. Therefore,
    $Y_G(\Zpos)$ is Zariski dense.
\end{proposition}

\begin{proof}
    For each $k \geq 1$, Lemma~\ref{lem:nonmonic-fiber-points}
    gives a positive integral point with
    $z = z_k = 1 + ab + a(ab)^{k - 1}c$.
    We check the discriminants and apply Lemma~\ref{lem:pell-equation}
    to obtain infinitely many positive integral points on
    infinitely many distinct fibers.

    Since $ab > 1$, the sequence $(z_k)$ is strictly increasing,
    and $z_k^2 - 4ab > 0$ for all sufficiently large $k$.
    Only finitely many of these discriminants are squares:
    if $z_k^2 - 4ab = w^2$, then
    $(z_k - w)(z_k + w) = 4ab$, which bounds $z_k$.
    Consequently, Lemma~\ref{lem:pell-equation} gives infinitely
    many positive integral points on $ax^2 - z_kxy + by^2 = -c$
    for all sufficiently large $k$.

    The surface $Y_G$ is geometrically integral, and the conics
    $ax^2 - z_kxy + by^2 = -c$ are smooth and irreducible for
    these values of $k$. A proper closed subset of an integral
    surface has only finitely many irreducible curve components.
    It cannot contain infinitely many distinct such conics,
    each with infinitely many points. Hence $Y_G(\Zpos)$ is
    Zariski dense.
\end{proof}

\section{Examples and applications}\label{sec:examples}

We apply the preceding results to the coefficient-free family
$xyz=x^m+y^n+1$ and to representation problems for binary quadratic
forms.

\subsection{The coefficient-free family}
\label{sec:extremal-family}
For $m,n\in\Zpos$, set
\[
    Y_{m,n}=(xyz=x^m+y^n+1),\qquad
    A(T)=T^m+1,\quad B(T)=T^n+1.
\]
Lemma~\ref{lem:2-integral-forms} identifies $\Z[Y_{m,n}]$ with
the rank-two cluster algebra with trivial coefficients. With
$x_0=x$ and $x_1=y$, its cluster variables satisfy
\begin{equation}\label{eq:rec-def-cv}
    x_{k-1}x_{k+1}=\begin{cases}
        1+x_k^m,&k\text{ even},\\
        1+x_k^n,&k\text{ odd}.
    \end{cases}
\end{equation}
Since $A^\dagger=A$ and $B^\dagger=B$, the groupoid
$\calg_{A,B}$ has one object. Write
\[
    \Gamma_{m,n}=\Gamma_{A,B}=\langle\mu_1,\mu_2\rangle,
    \qquad \mu_i=\mu_{A,B,i}.
\]

Theorem~\ref{thm:trichotomy} gives
\[
    \#Y_{m,n}(\Zpos)<\infty\quad\Longleftrightarrow\quad mn\leq3,
\]
recovering~\cite[Corollary A]{zhang-dynkin} when $n = 1$. The cardinalities for
$mn=1,2,3$ are respectively $5,6,9$.

\begin{proposition}\label{prop:extremal-affine-one-orbit}
    If $mn=4$, then $Y_{m,n}(\Zpos)$ is a single
    $\Gamma_{m,n}$-orbit.
\end{proposition}
\begin{proof}
    For $(m,n)=(2,2)$, the minimum argument in
    Proposition~\ref{prop:quadratic-affine-orbits} gives
    \[
        (\delta_-+\delta_+)u+\delta_-\delta_+=1.
    \]
    Thus $u=1$ and $\{\delta_-,\delta_+\}=\{0,1\}$, so every
    orbit contains $(1,1,3)$.

    For $(m,n)=(4,1)$, use the notation of
    Proposition~\ref{prop:quartic-affine-orbits}. Equation
    \eqref{eq:quartic-minimum-polynomial} becomes
    \[
        r(x^2-1)+tx=2x.
    \]
    If $x\geq2$, then $r=1$, $t=0$, and $(x-1)^2=2$, a
    contradiction. Hence $x=1$, $t=2$, and
    $\{\delta_-,\delta_+\}=\{1,2\}$. Every orbit therefore
    contains $(1,1,3)$ or $(1,2,2)$, which are exchanged by
    $\mu_2$. The case $(1,4)$ follows by interchanging $x$ and $y$.
\end{proof}

The unique orbit lies on $V(z-3)$ in $Y_{2,2}$ and on
$V(\Omega-5)$ in $Y_{4,1}$, where
$\Omega=yz^2-x^3z-x^2$. Neither $\mu_1$ nor $\mu_2$ has a
positive integral fixed point on these surfaces, so
$e_{\calo}=2$. Indeed, in bidegree $(2,2)$ a fixed point would
satisfy $x^2=y^2+1$ or $y^2=x^2+1$. In bidegree $(4,1)$,
$\mu_2$ would require $y^2=x^4+1$, whereas a fixed point of
$\mu_1$ would give $y=x^2-1$ and $y\mid2$, which is impossible.
Writing $E_{m,n}=E_{x^m+y^n+1}$, Propositions
\ref{prop:orbit-count-2-2} and~\ref{prop:orbit-count-4-1} give
\[
    E_{2,2}(N)=\frac{\log N}{\log\varphi}+O(1),
    \qquad \varphi=\frac{1+\sqrt5}{2},
\]
and
\[
    E_{4,1}(N)=E_{1,4}(N)
    =\frac{2\log N}{\log((5+\sqrt{21})/2)}+O(1).
\]
For nonreciprocal quartics, the corresponding orbit and counting
statements are those of \S\ref{sec:quartic-case}, with
$\Gamma_{A,B}$ in place of $\Gamma_{m,n}$ and $\varepsilon_A=1$.

\begin{proposition}\label{prop:extremal-positive-automorphisms}
    For $m,n\in\Zpos$, the group $\Aut_\Z^+(Y_{m,n})$ is
    generated by $\Gamma_{m,n}$ and, when $m=n$, the transposition
    $(x,y,z)\mapsto(y,x,z)$. In particular,
    \[
        [\Aut_\Z^+(Y_{m,n}):\Gamma_{m,n}]\leq2.
    \]
\end{proposition}
\begin{proof}
    By Lemma~\ref{lem:2-integral-forms}, the coordinate ring has
    the presentation
    \[
        \Z[u,u',v,v']/(uu'-v^n-1,\ vv'-u^m-1).
    \]
    Blanc--Dolgachev~\cite[Theorem~1]{blanc-dolgachev}, applied
    over $\Q$, express every automorphism as a diagonal
    automorphism followed by an element generated by
    $\Gamma_{m,n}$ and transposition when $m=n$. The diagonal
    factors on $u,v$ are rational roots of unity, hence signs.
    The remaining generators preserve positive integral points.
    Since this locus is nonempty, its preservation forces both
    signs to be $1$. All the remaining generators are defined
    over $\Z$.
\end{proof}

\begin{proposition}\label{prop:indefinite-3-2-orbits}
    The set $Y_{3,2}(\Zpos)$ contains infinitely many
    $\Gamma_{3,2}$-orbits. More precisely, for every integer
    $t\geq2$, the point
    \[
        P_t=(t^2+1,t^2+t+1,t^2-t+3)
    \]
    belongs to $Y_{3,2}(\Zpos)$ and is the unique point of
    minimum height in its orbit. The same conclusion holds for
    $Y_{2,3}$, and both surfaces have infinitely many
    $\Aut_\Z^+$-orbits of positive integral points.
\end{proposition}
\begin{proof}
    Write $P_t=(x_t,y_t,z_t)$. Direct expansion gives
    $x_ty_tz_t=x_t^3+y_t^2+1$. On the first two coordinates,
    \[
        \mu_1(x,y)=\left(\frac{y^2+1}{x},y\right),\qquad
        \mu_2(x,y)=\left(x,\frac{x^3+1}{y}\right).
    \]
    At $P_t$, the new coordinates are
    \[
        \frac{y_t^2+1}{x_t}=t^2+2t+2>y_t,\qquad
        \frac{x_t^3+1}{y_t}=(t^2+2)(t^2-t+1)>y_t.
    \]
    Each first mutation therefore increases the height. Moreover,
    \[
        x\leq y\ \Longrightarrow\ \frac{y^2+1}{x}>y,
        \qquad
        y\leq x\ \Longrightarrow\ \frac{x^3+1}{y}>x.
    \]
    Every reduced word alternates the two involutions. After
    its first mutation at $P_t$, each successive mutation
    replaces the smaller coordinate and strictly increases the
    height. Thus $P_t$ is the unique minimum in its orbit.
    Since $h(P_t)=t^2+t+1$ is strictly increasing, these orbits
    are distinct. Interchanging $x,y$ gives the result for
    $Y_{2,3}$; Proposition~\ref{prop:extremal-positive-automorphisms}
    identifies the full positive automorphism groups with the
    corresponding groups $\Gamma_{m,n}$.
\end{proof}

This proves the assertion about infinitely many automorphism
orbits in Conjecture~\ref{conj:indefinite-orbits} for these two
coefficient-free surfaces. The orbit question remains open in
bidegree $(3,2)$ for arbitrary positive coefficients.

\subsection{Mordell's example}
Recall the definition of $\calb_{m,n}$ from Section \ref{sec:ca}. Put 
\[
    \mathcal{E}_{m,n} = V(a_1,\ldots, a_{m-1}, b_1,\ldots, b_{n-1}).
\]
The ring of regular functions $\C[\mathcal{M}^{\rm ex}_{m,n}]$ of $\mathcal{M}^{\rm ex}_{m,n} := \pi^{-1}(\mathcal{E}_{m,n})$
is isomorphic to 
\[
    \C[x,y,z,a, b, c]/(xyz - a x^m - b y^n - c).
\]

Set $\P = {\rm Trop}(a,b,c)$ and $\Sigma_{m,n}^{\rm ex} = ({\bf x}, ({\bf p}_1,{\bf p}_2), M)$ where ${\bf p}_1 \in (\P\cup 0)^{n+1}$, ${\bf p}_2 \in (\P\cup 0)^{m+1}$ are
\[  
    {\bf p}_1 = (b,0,\ldots, 0,c), \qquad {\bf p}_2 = (c,0,\ldots, 0, a), \qquad M = \begin{pmatrix}
        0 & m \\ -n & 0
    \end{pmatrix}.
\]
The seed pattern $[\Sigma_{m,n}^{\rm ex}] = \{\Sigma_k = ({\bf x}_k, ({\bf p}_{k;1}, {\bf p}_{k;2}), M_k) : k \in \Z\}$ is a {\it seed pattern of geometric type} in the sense of \cite{fomin-zelevinsky-1}. Note here that we allow some coefficients to be equal to zero: see \cite{nakanishi:addendum} for this subtlety. Indeed, \eqref{eq:coef-mut} guarantees that 
\[
    p_{k;i,s} = 0 , \qquad \forall k \in \Z, i \in \{1,2\}, s \in [1,d_i-1], 
\]
where $d_1 = n$ and $d_2 = m$, so the normalisation condition \eqref{eq:normalised-coeffs} becomes
\begin{equation}\label{eq:geometric-type}
    p_{k;i}^- \oplus p_{k;i}^+ = 1, \qquad \forall k\in \Z, i \in \{1,2\}.
\end{equation}
Recall the notation from \S \ref{sec:appendix-rank2}. In particular, $\{x_k : k \in \Z\}$ denotes the set of cluster variables. The exchange relations are precisely 
\begin{equation}\label{eq:exchange-mordell}
    x_{k} \, x_{k+2} = p_{k+1}^+ x_{k+1}^{d_{k+1}} + p_{k+1}^-,
\end{equation}
It follows from \cite{fomin-zelevinsky-3} that 
\[
    \C[\mathcal{M}^{\rm ex}_{m,n}] \cong \mathcal{U}_\LL(\Sigma_{m,n}^{\rm ex}) \otimes \C,
\]
where $\mathcal{U}_\LL(\Sigma_{m,n}^{\rm ex})$ is the upper cluster algebra over the ground ring $\LL = \Z[a,b,c]$ associated to $[\Sigma_{m,n}^{\rm ex}]$.

It is well-known that the coefficients $p_{k;i}^\pm$ of $[\Sigma_{m,n}^{\rm ex}]$ can be conveniently read-off from the {\it extended exchange matrix}, as we now explain. 
For each $k \in \Z, i \in \{1,2\}$ there exist unique $P_{k;i}^{\pm} \in \Z^3$ such that 
\[
    p_{k;i}^\pm = (a, b, c)^{P_{k;i}^\pm }.
 \]
 The normalisation condition \eqref{eq:geometric-type} implies $P_{k;i}^{\pm}\in (\Znn)^3$ and that if
 \[
    P_{k;i}^\pm = \begin{pmatrix}
        \lambda_{k;i}^\pm \\ \mu_{k;i}^{\pm} \\ \nu_{k;i}^\pm
    \end{pmatrix},
 \]
 then 
\begin{equation}\label{eq:normalisation-separation}
\lambda_{k;i}^+ \lambda_{k;i}^- = \mu_{k;i}^+\mu_{k;i}^-  =  \nu_{k;i}^+ \nu_{k;i}^- = 0.
\end{equation}
Define the so-called extended exchange matrix $\widetilde{M}_k \in M_{5,2}(\Z)$ by 
\[
   \widetilde{M}_k = \begin{pmatrix}
       M_k \\ C_k
   \end{pmatrix} 
\]
where if $C_{k;i}$ denotes the $i^{\rm th}$ column of $C_k$, then $C_{k;i} = P_{k;i}^+ - P_{k;i}^-$. That $C_k$ completely determines $P_{k;i}^\pm$ is guaranteed by \eqref{eq:normalisation-separation}. 

Recall that for a $(2+d)\times 2$ matrix $\widetilde{M} = (m_{i,j})$ whose principal $2\times 2$ submatrix $M$ is a matrix of the form described in Definition \ref{def:seed-rank2} and $k \in \{1,2\}$, the {\it mutation of $\widetilde{M}$ in direction $k$} is $\mu_k(\widetilde{M}) = (m_{i,j}')$, where
\begin{equation}\label{eq:matrix-mut}
m_{i,j}' = \begin{cases}
            -m_{i,j}, & \text{ if } i = k \text{ or } j = k, \\
            m_{i,j} + [m_{i,k}]_+[m_{k,j}]_+ - [-m_{i,k}]_+[-m_{k,j}]_+, &\text{ otherwise},
        \end{cases}
\end{equation}
with $[x]_+ := \max(0,x)$. Note that $\mu_k^2(\widetilde{M}) = \widetilde{M}$. It was shown in \cite{fomin-zelevinsky-1} that $\widetilde{M}_{k'} = \mu_i(\widetilde{M}_k)$ whenever {\begin{xy}(0,1)*+{k}="A",(10,1)*+{k'}="B",\ar@{-}^i"A";"B" \end{xy}}. Concretely, we have 
\[
    \widetilde{M}_0 = \begin{pmatrix}
        0 & m \\ -n & 0 \\ 0 & 1 \\ -1 & 0 \\ 1 & -1
    \end{pmatrix}, \quad \widetilde{M}_1 = \begin{pmatrix}
        0 & -m \\ n & 0 \\ 0 & 1 \\ 1 & 0 \\ -1 & m-1
    \end{pmatrix}, \quad \widetilde{M}_2 = \begin{pmatrix}
        0 & m \\ -n & 0 \\ n & -1 \\ 1 & 0 \\ nm-n-1 & 1-m
    \end{pmatrix}, \ldots
\]

The general formula is given by 
\[
    C_{2k+2;1} = n[C_{2k+1;2}]_+ - C_{2k;1}  \quad \text{and}\quad C_{2k+1;2} = m[C_{2k;1}]_+ - C_{2k-1;2}.
\]
\begin{example}\label{eq:Mordell-mn<4}
    One checks by a direct calculation (e.g. using \cite{visualca}) that $\widetilde{M}_{10} = \widetilde{M}_0$ when $mn =1$, $\widetilde{M}_6 = \widetilde{M}_0$ when $mn =2$, and $\widetilde{M}_8 = \widetilde{M}_0$ when $mn =3$. 
\end{example}
\begin{example}
    Consider the case $(m,n) = (2,2)$. For $\ell \geq 1$ the extended exchange matrix $\widetilde{M}_{\ell}$ is given by
    \[
     \begin{pmatrix}
        0 & 2 \\ -2 & 0\\
        \ell & -(\ell-1)\\
        \ell-1 & - (\ell -2)\\
        1 & -1
    \end{pmatrix} \text{ if } \ell \text{ even},  \quad \begin{pmatrix}
        0 & -2 \\ 2 & 0\\
        -(\ell-1) &\ell  \\
        - (\ell -2) & \ell-1 \\
        -1 & 1
    \end{pmatrix} \text{ if } \ell \text{ odd}.
    \]
    The exchange relations \eqref{eq:exchange-mordell} are then 
    \[
        x_0 x_2 = bx_1^2 + c,
    \]
    and for $k \geq 1$,
    \[
        x_k x_{k+2} = a^kb^{k-1}c + x_{k+1}^2.
    \]
\end{example}

\begin{example}
    When $(m,n) = (3,3)$, the matrices $\widetilde{M}_1, \widetilde{M}_2, \widetilde{M}_3$ and $\widetilde{M}_4$ are respectively 
    \[
     \begin{pmatrix}0&-3\\3&0\\0&1\\1&0\\-1&2\end{pmatrix}, \quad \begin{pmatrix}0&3\\-3&0\\3&-1\\1&0\\5&-2\end{pmatrix}, \quad \begin{pmatrix}0&-3\\3&0\\-3&8\\-1&3\\-5&13\end{pmatrix}, \quad \begin{pmatrix}0&3\\-3&0\\21&-8\\8&-3\\34&-13\end{pmatrix}.
    \]
    The column vectors $C_{0;1}, C_{1;2}, C_{2;1}, C_{3;2}, C_{4;1}$, etc. are precisely the vectors of exponents computed in \cite[(14)]{mordell-1952}, and the exchange relations \eqref{eq:exchange-mordell} are given in \cite[(15)]{mordell-1952}.
\end{example}
It is clear from Example \ref{eq:Mordell-mn<4} that Mordell's method can not possibly furnish infinitely many integer solutions when $mn \leq 3$. 

\subsection{Representation problems}\label{sec:rep-problems}
For a pure polynomial $G\in \Znn[x,y]$ with $c>0$, set
\[
    F_z(x,y)=G(x,y)-c-zxy=A(x)+B(y)-2c-zxy.
\]
The positive integral points of $Y_G$ with third coordinate $z$
are the positive integral representations $F_z(x,y)=-c$.
In bidegree $(2,2)$, $F_z$ is a quadratic polynomial. It is
homogeneous precisely when $G=ax^2+by^2+c$, in which case
\[
    F_z(x,y)=ax^2-zxy+by^2.
\]

\begin{proposition}\label{prop:indefinite-form}
    Let $G=ax^2+by^2+c$, with $a,b,c\in\Zpos$.
    For every $(x,y,z)\in Y_G(\Zpos)$, the binary quadratic form
    $F_z$ is indefinite.
\end{proposition}
\begin{proof}
    The surface equation gives
    \[
        z=\frac{ax^2+by^2+c}{xy}
        >a\frac{x}{y}+b\frac{y}{x}\geq2\sqrt{ab},
    \]
    so the discriminant $D=z^2-4ab$ is positive.
\end{proof}

Retain $G=ax^2+by^2+c$. The isomorphisms of
Example~\ref{ex:2-2-isoms} preserve $z$. If $F_z(x,y)=-c$, then
$u=yz-ax$ and $v=xz-by$ satisfy
\[
    u^2-zuy+aby^2=-ac,\qquad
    abx^2-zxv+v^2=-bc.
\]
Starting from $(x,y,z)=(1,1,a+b+c)$ and iterating these maps
shows that the fixed form $x^2-zxy+aby^2$, with $z=a+b+c$,
represents every integer $-a(ab)^{k-1}c$, $k\geq1$.

The normalisation in Theorem~\ref{thm:affine-counting} is $a=b=1$. In this case, Theorem~\ref{thm:trichotomy}(ii), Theorem~\ref{thm:affine-counting}, and \S\ref{sec:quadratic-case} give only finitely many values of $z$,
and each $\Gamma_{A,B}$-orbit has fixed discriminant $z^2-4$
and logarithmic growth. If $ab>1$, Proposition
\ref{prop:nonmonic-quadratic-density} instead gives infinitely
many values of $z$ for which $F_z$ represents $-c$ in infinitely
many positive ways, and $Y_G(\Zpos)$ is Zariski dense.

For $ab>1$, Lemma~\ref{lem:nonmonic-fiber-points} gives the
sequence $z_k=1+ab+a(ab)^{k-1}c$. One may ask how the ideal class
groups of the quadratic orders
$\Z[\sqrt{z_k^2-4ab}]$ vary for sufficiently large $k$.
Another question is whether the infinitely many automorphism
orbits predicted by Conjecture~\ref{conj:indefinite-orbits} admit
an interpretation in terms of representations by binary
polynomials.

For $Y_c=(xyz=x^2+y^2+c)$, the fibers of $z$ are norm curves.
This identifies the group action underlying the Pell equations
in Lemma~\ref{lem:pell-equation}.

\begin{proposition}\label{prop:norm-torus-curves}
    Let $K$ be a field of characteristic different from $2$,
    let $c\in K^\times$, and fix $z_0\in K$ with
    $D=z_0^2-4\ne0$. Set $L=K[t]/(t^2-D)$.
    Write $\operatorname{Res}_{L/K}\G_m$ for the restriction
    of scalars of $\G_m$, whose $K$-points are $L^\times$,
    and $\rmn_{L/K}$ for the norm. Set
    \[
        T=\ker\left(\rmn_{L/K}:\operatorname{Res}_{L/K}\G_m
        \longrightarrow\G_m\right).
    \]
    The curve $Y_c\cap V(z-z_0)$ is a torsor under the
    one-dimensional torus $T$. A $K$-point identifies the curve
    with $T$, and the curve becomes isomorphic to $\G_m$ over
    $K(\sqrt D)$.
\end{proposition}
\begin{proof}
    The algebra $L$ is quadratic \'etale. Let $\delta$ be the
    image of $t$ in $L$ and set $\eta=(z_0+\delta)/2$. Since
    $\eta+\overline\eta=z_0$ and $\eta\overline\eta=1$,
    \[
        \rmn_{L/K}(x-\eta y)=x^2-z_0xy+y^2.
    \]
    The basis $1,\eta$ of $L$ thus identifies the curve with
    $\rmn_{L/K}(u)=-c$. Multiplication by $T$ makes this a
    torsor: any two solutions differ by a unique norm-one
    factor. A solution $u_0\in L$ identifies $T$ with the
    curve by $v\mapsto u_0v$.

    Over $K(\sqrt D)$, the norm equation becomes $rs=-c$,
    which is isomorphic to $\G_m$ via $r\mapsto(r,-c/r)$.
\end{proof}

\begin{remark}\label{rem:fibonacci-fibers}
    For $K=\Q$ and $z_0=3$, the fibers with $c=1$ and $c=-1$
    are respectively the norm $-1$ and norm $1$ curves in
    $\Q(\sqrt5)$. Multiplication by the norm-one unit
    $(3+\sqrt5)/2$ induces $(x,y)\mapsto(y,3y-x)$, so successive
    coordinates satisfy $u_{j+2}=3u_{j+1}-u_j$.
    Starting from $(1,2)$ for $c=1$ gives
    $1,2,5,13,\ldots$; starting from $(1,3)$ for $c=-1$ gives
    $1,3,8,21,\ldots$. These are respectively the odd-indexed
    and even-indexed Fibonacci numbers.
\end{remark}

\appendix
\section{Background on generalised cluster algebras}
\label{sec:appendix}
\subsection{Background}\label{ss:appendix-background}

We recall standard definitions and basic properties of generalised cluster algebras, following closely \cite{burcroff-Lee-mou:generalised}. 

Fix $\P$ a torsion-free multiplicative abelian
group, which we call the {\it coefficient group}, fix $d$ a positive integer, write $I = \{1,\ldots, d\}$ and let $\F$ be a field isomorphic to the field of rational functions in $d$ independent variables over the fraction field $\Q\P$ of the group ring $\Z\P$. 

\begin{definition}\label{def:seed-rank2}
A \emph{seed} in $\F$ is a triple
$\Sigma=({\bf x},{\bf p},M)$ satisfying the following conditions.
\begin{enumerate}
    \item The {\it cluster} ${\bf x}=(x_1,\ldots , x_d)$ is algebraically independent over
    $\Q\P$ and $\F=\Q\P(x_1,\ldots, x_d)$.
    \item The {\it exchange matrix} $M=(m_{i,j})$ is a skew-symmetrisable\footnote{that is, if there exists a diagonal matrix $D$ with positive integer diagonal entries such that $DM$ is skew-symmetric} $d \times d$ integer matrix.
    \item The {\it coefficient tuples} ${\bf p} = ({\bf p}_1, \ldots, {\bf p}_d)$ consist, for each $i \in I$, of a $(d_i+1)$-tuple ${\bf p}_i=(p_{i,0},\ldots,p_{i,d_i})\in\P^{d_i+1}$, where $d_j = \gcd(m_{k,j} : k \in I)$ for a nonzero column
    and $d_j = 1$ for a zero column.
\end{enumerate}
\end{definition}
By the {\it rank} of a seed $\Sigma=({\bf x},{\bf p},M)$ in $\F$ we mean the length $d$ of the cluster ${\bf x}$. Elements $x \in {\bf x}$ are called {\it cluster variables} and elements $p \in {\bf p}_i, i \in I$ {\it coefficients} of $\Sigma$. We write $p_i^-=p_{i,0}$ and $p_i^+=p_{i,d_i}$.

For $k \in I$, define $x_k' \in \F$ via the {\it exchange relation}
\begin{equation}\label{eq:generalised-exchange}
    x_k\, x_k'  = \sum_{r = 0}^{d_k} p_{k,r} \;\prod_{j \in I} x_j^{\lambda(r,j)} ,
\end{equation}
where 
\[
    \lambda(r,j) = \begin{cases}
        \frac{r}{d_k}m_{j,k}, & \text{if }m_{j,k} \geq 0, \\
        \frac{r-d_k}{d_k}m_{j,k}, & \text{if }m_{j,k} < 0.
    \end{cases} 
\]

\begin{definition}\label{def:mutation}
Consider two seeds $\Sigma=({\bf x},{\bf p},M)$ and $\Sigma'=({\bf x}',{\bf p}',M')$ in $\F$. We say that $\Sigma'$ is obtained from $\Sigma$ via {\it mutation} in direction
$k \in I$ if the following conditions hold.

1) ${\bf x}'$ is obtained from ${\bf x}$ by replacing $x_k \in {\bf x}$ with $x_k'$ as defined in \eqref{eq:generalised-exchange}. 

2) The coefficients satisfy
$p'_{k,r}=p_{k,d_k-r}$ and, for $i \in I, i\neq k$,
\begin{equation}\label{eq:coef-mut}
    \frac{p'_{i,r}}{p'_{i,0}}
    =\frac{p_{i,r}}{p_{i,0}}
      \left(\frac{p_k^+}{p_k^-}\right)^{r[m_{k,i}]_+/d_i}
      (p_k^-)^{rm_{k,i}/d_i}
      \qquad(0\leq r\leq d_i).
\end{equation}

3) $M' = \mu_k (M)$, where $\mu_k$ denotes matrix mutation in direction $k$ (see e.g. \cite[(7.5)]{burcroff-Lee-mou:generalised}).
\end{definition}

Note that condition 2) of Definition \ref{def:mutation} does not uniquely specify (the coefficient tuples of) $\Sigma'$ in terms of (the coefficient tuples of) $\Sigma$. We shall return to this in the next section.

Consider $\T_d$ the $d$-regular tree, where edges incident to every vertex are labeled by $I$. We write \tkt to mean that the two vertices $t$ and $t'$ of $\T_d$ are connected by an edge labeled by $k \in I$.

A {\it seed pattern} of rank $d$ in $\F$ is an assignment 
\[
\{\Sigma_t : t \in \T_d \}, \quad \Sigma_t = ({\bf x}_t, {\bf p}_t, M_t)
\]
of a seed $\Sigma_t$ of rank $d$ in $\F$ to each $t \in \T_d$ such that $\Sigma_{t'}$ is obtained from $\Sigma_{t}$ by mutation in direction $k$ whenever \tkt. For each $t \in \T_d, i \in I$ and $0 \leq r \leq d_i$, the corresponding cluster variable and coefficient of $\Sigma_t$ are denoted $x_{t;i}$ and $p_{t;i,r}$ respectively. 

In applications, it is often useful to work over rings that are smaller than $\Z\P$. By a {\it ground ring} we mean a subring $\LL \subset \Z\P$ containing $\{p_{t;i,r} : t \in \T_d, i \in I, 0 \leq r \leq d_i\}$.

\begin{definition}\label{def:gca}
    The {\it generalised upper cluster algebra} over the ground ring $\LL$ and associated to the seed pattern $[\Sigma] = \{\Sigma_t : t \in \T_d\}$ in $\F$ is 
    \[
        \mathcal{U}_\LL(\Sigma) = \bigcap_{t \in \T_d} \LL[x_{t;1}^{\pm 1}, \ldots,x_{t;d}^{\pm 1}].
    \]
    When $\LL = \Z\P$, we simply write $\mathcal{U}(\Sigma)$.
\end{definition}
We refer to $[\Sigma]$ as the seed pattern underlying $\calu_\LL(\Sigma)$, and we call a {\it cluster} (respectively, a {\it cluster variable}, a {\it coefficient tuple} or a {\it coefficient}) of $\calu_\LL(\Sigma)$ a cluster (respectively, a cluster variable, a coefficient tuple or a coefficient) of some seed in $[\Sigma]$. We write 
\[
    \mathfrak{X}([\Sigma]) = \{x_{t;i} : t \in \T_d , i\in I\},
\]
for the set of cluster variables of $[\Sigma]$. A generalised upper cluster algebra and its underlying seed pattern are said to be of {\it finite cluster-type} if their set of cluster variables is finite. These were classified by Fomin and Zelevinsky in \cite{fomin-zelevinsky-2}, in terms of root systems. 

The following remarkable statement, known as the {\it positive Laurent phenomenon}, was conjectured in \cite{fomin-zelevinsky-1}, and proved in increasing levels of generality in \cite{chekhov-shapiro, GHKK,burcroff-Lee-mou:generalised}.

\begin{theorem}\cite[Theorem 7.9]{burcroff-Lee-mou:generalised}\label{thm:pos-laurent-phen}
    Let $[\Sigma]$ be a seed pattern of rank $d$. For each $x \in  \mathfrak{X}([\Sigma])$, we have 
    \[
         x \in  \bigcap_{t \in \T_d} \Znn\P[x_{t;1}^{\pm 1}, \ldots, x_{t;d}^{\pm 1}].
    \]
\end{theorem}

\subsection{Normalisations}
Keep the notation of \S \ref{ss:appendix-background}, and equip $\P$ with an {\it auxiliary addition} $\oplus$ which is commutative, associative, and distributes over the multiplication. We say that a seed $\Sigma$ (resp. a seed pattern) in $\F$ is {\it normalised} if the coefficient tuples ${\bf p}_i, i \in I$ of $\Sigma$ (resp. of each of its seeds $\Sigma$) satisfy
\begin{equation}\label{eq:normalised-coeffs}
     \bigoplus_{r=0}^{d_i}p_{i,r}=1
        \qquad \forall i \in I.
    \end{equation}
Then, mutation of seeds is deterministic: the ratios in \eqref{eq:coef-mut} together with \eqref{eq:normalised-coeffs} determine ${\bf p}'_i$ uniquely. In other words, a normalised seed pattern is uniquely determined by a single normalised seed (see \cite[Proposition 5.4]{fomin-zelevinsky-1}). Fix $t_0 \in \T_d$, which we call the root, and a normalised seed $\Sigma$ in $\F$. Then,  $[\Sigma] = \{\Sigma_t : t \in \T_d\}$ denotes the unique normalised seed pattern in $\F$ such that $\Sigma_{t_0} = \Sigma$; $\Sigma$ is called the {\it initial seed} of $[\Sigma]$.

Now let $\P$ be the free abelian multiplicative group with basis $S=\{f_1,\ldots, f_q\}$, whose elements we call {\it frozen variables} of $\P$. Define $\oplus$ via
\begin{equation}\label{eq:trop-sum}
    \paren{\prod_{j=1}^q f_j^{a_j}} \oplus \paren{\prod_{j=1}^q f_j^{b_j}} = \prod_{j=1}^q f_j^{\min(a_j,b_j)}.
\end{equation}
The pair $(\P,\oplus)$ is called {\it a tropical semifield}, denoted ${\rm Trop}(S)$. Here \eqref{eq:normalised-coeffs} is equivalent to saying that, for each $i$, the $p_{i,r}$ are monomials (as opposed to a Laurent monomial) in $f_1,\ldots, f_q$ with no common factor. 

When the underlying coefficient group $\P$ is a tropical semifield, normalised seeds, normalised seed patterns, and their corresponding generalised upper cluster algebras are said to be {\it of geometric type}. The {\it frozen variables} of $[\Sigma]$ are the frozen variables of $\P$; we denote this set by ${\rm Froz}(\Sigma)$. There are two common choices of ground rings: 
\[
    \LL = \Z\P \quad \text{and} \quad \LL = \Z[f_1,\ldots , f_q].
\]
In the former (resp. the latter), $\calu_\LL(\Sigma)$ is said to have {\it inverted} (resp. {\it non-inverted}) frozen variables.

The following strengthening of Theorem \ref{thm:pos-laurent-phen} is sometimes referred to as the {\it strong positive Laurent phenomenon}, see \cite[\S 11]{fomin-zelevinsky-2}.

\begin{theorem}\label{thm:appendix-strong-pos-laur-phen}
    Let $[\Sigma]$ be a seed pattern of rank $d$ and of geometric type. For each $x \in \mathfrak{X}([\Sigma])$, we have 
    \[
        x \in \bigcap_{t \in \T_d} \Znn[x_{t;1}^{\pm 1}, \ldots, x_{t;d}^{\pm 1}, f_1,\ldots, f_q].
    \]
\end{theorem}
In other words, each cluster variable of $[\Sigma]$ is a Laurent polynomial in the variables of a cluster, {\it polynomial} in the frozen variables, and with non-negative coefficients. 

Before proving Theorem \ref{thm:appendix-strong-pos-laur-phen}, we introduce a bit of notation. The {\it ambient semifield} of $[\Sigma]$, denoted $\F_{>0}$, is the smallest semifield in $\F$ containing $\Q_{>0}$, $\mathfrak{X}([\Sigma])$ and ${\rm Froz}(\Sigma)$; it is clear from Equation \ref{eq:generalised-exchange} that 
\[
    \F_{>0} = \Q_{>0}(f_1,\ldots, f_q, {\bf x}),
\]
where the latter is the set of all rational expressions in $f_1,\ldots, f_q$ and ${\bf x}$ with coefficients in $\Z_{\geq 0}$. 

\begin{lemma}\label{lem:sf-homo}
    For $[\Sigma]$ as in Theorem \ref{thm:appendix-strong-pos-laur-phen}, there exists a unique semifield homomorphism $\F_{>0} \longrightarrow \P$ which sends each cluster variable of $[\Sigma]$ to $1 \in \P$ and fixes $f_1,\ldots, f_q$.
\end{lemma} 
\begin{proof}
     To see this, fix a cluster ${\bf x}$, and define the semifield homomorphism $Z \mapsto Z\vert_{\P}$ by $x \mapsto 1$ for $x \in {\bf x}$ and by fixing ${\rm Froz}(\Sigma)$. Then, for every $k \in I$, Equation \ref{eq:generalised-exchange} and Equation \ref{eq:normalised-coeffs} imply 
    \begin{align*}
        x_k'\big\vert_{\P} &= \frac{\sum_{r = 0}^{d_k} p_{k,r} \;\prod_{j \in I} x_j^{\lambda(r,j)}}{x_k}\Big\vert_{\P}\\
        &= \bigoplus_{r=0}^{d_k} p_{k,r} = 1.
    \end{align*}
    By induction on the distance from $t_0$ in $\T_d$, the homomorphism $Z \mapsto Z\vert_\P$ is readily seen to satisfy the properties in the claim. Its uniqueness is obvious. 
\end{proof}

{\it Proof of Theorem \ref{thm:appendix-strong-pos-laur-phen}.} Fix $x \in \mathfrak{X}([\Sigma])$ and $t \in \T_d$. By Theorem \ref{thm:pos-laurent-phen},  $x \in \Znn\P[x_{t;1}^{\pm 1}, \ldots, x_{t;d}^{\pm 1}]$. Write
\begin{equation}\label{eq:Laurent-expansion-x}
    x = \frac{N_1(x_{t;1},\ldots, x_{t;d})}{N_2}, \quad N_2 = \prod_{i} x_{t;i}^{s_i} \prod_{j} f_j^{t_j} 
\end{equation}
for some unique $N_1,N_2$ coprime with $N_1 \in \Znn[{\bf x}, f_1,\ldots, f_q]$. Applying the homomorphism of Lemma \ref{lem:sf-homo} to $x$ gives $1 = \prod_{j} f_j^{-t_j}$, and so $t_1 = \ldots, t_q = 0$, which concludes the proof. \qed 

\subsection{The case of rank \texorpdfstring{$2$}{2}}
\label{sec:appendix-rank2}
Let $[\Sigma]$ be a seed pattern of rank $2$. 
We label the rooted 2-regular tree $(\T_2,t_0)$ as follows:
\begin{equation}\label{eq:right}
\begin{xy}
(-14,1)*+{\cdots}="A1",(0,1)*+{t_{-1}}="A2",(14,1)*+{t_{0}}="A3", (28,1)*+{t_{1}}="A4",(42,1)*+{t_2}="A5",
(58,1)*+{t_3}="A6",(74,1)*+{t_4}="A7", (88,1)*+{\cdots}="A8",
\ar@{-}^{1}"A1";"A2", \ar@{-}^{2}"A2";"A3", \ar@{-}^{1}"A3";"A4",
\ar@{-}^{2}"A4";"A5",\ar@{-}^{1}"A5";"A6",\ar@{-}^{2}"A6";"A7",\ar@{-}^{1}"A7";"A8",
\end{xy}
\end{equation}

Since the vertices of $\T_2$ are in bijection with $\Z$, we write $\Sigma_k := \Sigma_{t_k}$
and
\[
    {\bf x}_{t_k} = (x_{k;1}, x_{k;2}), \quad {\bf p}_{t_k;i} = (p_{k;i,r})_{r=0}^{d_i}, \quad M_k = (-1)^k M_0.
\]
Note that $x_{k;i} = x_{k-1;i}$ and $p_{k;i,r} = p_{k+1;i,d_i-r}$ whenever $k + i$ is odd, in which case  we further simplify notation and set 
\[
    x_{k} = x_{k;i} \quad \text{and} \quad {\bf p}_{k+1} = {\bf p}_{k+1;i}.
\]
Extend $d_k$ 2-periodically by setting $d_{k+2} = d_k$. Each ${\bf p}_{k}:= (p_{k;r})_{r=0}^{d_{k}}$ is a $(d_k+1)$-tuple.  
Assume that $M_0 \neq 0$;
after interchanging the two indices if necessary,
we may suppose that
\[
  M_0 = \begin{pmatrix}
      0 & d_2 \\ -d_1 & 0
  \end{pmatrix},  
\]
for some $d_1, d_2 >0$. The exchange relations are then precisely
\[
    x_{k} x_{k+2} = E_{k+1}(x_{k+1}) = \sum_{r=0}^{d_{k+1}} p_{k+1;r} (x_{k+1})^r.
\]
Each $E_k$ appearing in an exchange relation is called an {\it exchange polynomial}. 

If $[\Sigma]$ is of geometric type (i.e. $\P$ is a tropical semifield and $[\Sigma]$ is normalised), the coefficients are determined recursively by
\begin{equation}\label{eq:appendix-rank-two-coefficients}
 p_{k+2;r} = 
 \frac{p_{k;d_k-r}\,(p_{k+1}^-)^{d_k-r}}
 {\displaystyle\bigoplus_{\ell=0}^{d_k}p_{k;\ell}(p_{k+1}^-)^{\ell}}
 \qquad(0\leq r\leq d_k),
\end{equation}
with $p_{k+1}^- := p_{k+1;0}$. 

We will make use of the following example in \S \ref{sec:gca}.
\begin{example}\label{ex:non-normalised-greedy}
Consider a seed pattern $[\Sigma] = \{\Sigma_t : t \in \T_2\}$ of rank $2$, such that $\Sigma_{t_0} = ({\bf x}, {\bf p}, M)$ (for some $t_0 \in \T_2$) where 
\[
    {\bf x} = (x,y), \quad {\bf p} = ((p_{1,0}, \ldots p_{1,n-1}, 1),(1, p_{2,1},\ldots, p_{2,m})), \quad M = \begin{pmatrix}
        0 & m \\ -n & 0
    \end{pmatrix}.
\]
Let $\{x_k : k \in \Z\}$ denote the set of cluster variables, and write $x = x_0$, $x' = x_2$, $y = x_1$ and $y' = x_{-1}$. We have the exchange relations 
\begin{equation}\label{eq:appendix-one-step}
        x x' = 1 + \sum_{r=1}^n p_{1,n-r}y^r\quad \text{and} \quad y y' = 1 + \sum_{r=1}^m p_{2,r}x^r.  
\end{equation}

For $(u,v) \in \Z^2$, let $x[u,v]$ be the corresponding greedy element, as defined in  \cite[Definition 3.10]{burcroff-Lee-mou:generalised}. We have
\[
 y'= x[0,1], \; x = x[-1,0], \quad y = x[0,-1], \; x' = x[1,0]. 
\]
Every cluster variable is, up to scaling by $\P$, a greedy
element. The formula for compatible gradings and its behaviour
under mutation
\cite[Theorem~3.18 and Proposition~A.9]
{burcroff-Lee-mou:generalised}, together with rescaling
by elements of $\P$
\cite[Proposition~7.13 and Remark~7.15]
{burcroff-Lee-mou:generalised}, give
\[
    x[u,v] \in \bigcap_{k \in \Z} \Znn\P[x_k^{\pm 1}, x_{k+1}^{\pm 1}].
\]
Of particular interest in this article is the element $x[1,1]$. Note here that the maximal Dyck path has a single horizontal edge $h_1$ and a single vertical edge $v_1$. The compatibility condition implies $\omega(h_1)\omega(v_1) = 0$. The corresponding gradings contribute $1, p_{1,n-r}y^r $ for each $1 \leq r \leq n$ and $p_{2,s}x^s$ for each $1 \leq s \leq m$. Thus,
\begin{equation}
    x[1,1] = \frac{1 + \sum_{r=1}^m p_{2,r}x^r + \sum_{s=1}^n p_{1,n-s}y^s}{xy}.
\end{equation}    

If $\P = {\rm Trop}(a_1,\ldots, a_m,b_1,\ldots, b_n)$, and 
\[
    {\bf p} = ((b_n,\ldots, b_1,1),(1,a_1,\ldots, a_m)),
\]
then $[\Sigma]$ is a seed pattern of geometric type, and
\begin{equation}
    x[1,1] = \frac{1 + \sum_{i=1}^m a_i x^i + \sum_{j=1}^n b_j y^j}{xy}.
\end{equation}
\end{example}

\printbibliography

\end{document}